\documentclass[11pt]{amsart}
\usepackage{hyperref}
\hypersetup{nesting=true,debug=true,naturalnames=true}
\usepackage{graphicx,amssymb,upref}

\usepackage[usenames,dvipsnames]{pstricks}
\usepackage{epsfig}
\usepackage{color}
\usepackage[latin1]{inputenc}
\usepackage[T1]{fontenc}
\usepackage{amsmath, amsthm}
\usepackage[english]{babel}
\usepackage{amssymb}
\usepackage{latexsym}
\usepackage{graphicx,amssymb,upref}
\usepackage[all]{xy}
\usepackage{hyperref}
\hypersetup{nesting=true,debug=true,naturalnames=true}
\usepackage{tikz}
\usetikzlibrary{matrix}

\newtheorem{theorem}{Theorem}[section]

\theoremstyle{definition}
\newtheorem{definition}[theorem]{Definition}
\newtheorem{example}[theorem]{Example}

\newcommand{\ot}{\otimes}
\newcommand{\co}{\circ}

\let\<\langle
\let\>\rangle

\let\uml\"

\title{Distributive laws and Hopf quasigroups} 

\begin{document}

\maketitle

\begin{center}
	{\bf Ram\'on
		Gonz\'{a}lez Rodr\'{\i}guez}
\end{center}

\vspace{0.1cm}
\begin{center}
	{\small \vspace{0.4cm}  [https://orcid.org/0000-0003-3061-6685]}
\end{center}
\begin{center}	{\small  CITMAga, 15782 Santiago de Compostela, Spain}
\end{center}
\begin{center}
	{\small  Universidade de Vigo, Departamento de Matem\'{a}tica Aplicada II,   E-36310 
		Vigo, Spain\\email: rgon@dma.uvigo.es}
\end{center}

\vspace{0.1cm}

\begin{abstract}  
 In this paper we introduce the notion of $a$-monoidal distributive law between two Hopf quasigroups $A$ and $H$. We prove that every $a$-monoidal distributive law induce a product on $A\ot H$, called the wreath product, thanks to which $A\ot H$ becomes in a Hopf quasigroup. Finally, using this construction, we show that double cross products of Hopf quasigroups, cross products of Hopf quasigroups with a skew pairing between them, Hopf quasigroups defined by the twisted double method, smash products of Hopf quasigroups and twisted smash products of Hopf quasigroups are examples of wreath products associated to $a$-monoidal distributive laws.
 \end{abstract} 

\vspace{0.2cm} 



\section{Introduction} The theory of distributive laws between monads was initiated by Beck \cite{Beck} and Barr \cite{Barr}  in
the seventies of the last century. A distributive law between two monoids $A$ and $H$ in a symmetric monoidal category ${\mathcal C}$ is a morphism $\Psi:H\ot A\rightarrow A\otimes H$ which is compatible with the monoid structures. It is well-known that a distributive law $\Psi:H\ot A\rightarrow A\otimes H$ induces  a monoid structure on the tensor product $A\otimes H$  commuting with  the action associated to the product of $A$ on the left and with  the action associated to the product of $H$ on the right. This monoid, denoted by  $A\ot_{\Psi}H$, admits several names in the literature (wreath product, twisted product, smash product, etc.) and its unit and product are defined by $$\eta_{A\ot_{\Psi}H}=\eta_{A}\ot \eta_{H},\;\;\;\mu_{A\ot_{\Psi}H}=(\mu_{A}\ot \mu_{H})\co (id_{A}\ot \Psi\ot id_{H}),$$ 
where $\eta_{A}$, $\eta_{H}$ denote the units of $A$ and $H$ and $\mu_{A}$, $\mu_{H}$ the corresponding products.

On the other hand,  bimonoids are monoids in the category of comonoids. From this point of view, a distributive law in the category of comonoids is a distributive law between the underlying monoids satisfying that  is a comonoid morphism. This kind of  distributive laws induce a wreath product bimonoid, where the product is the wreath product and the comonoid structure is the one associated to the tensor product comonoid. A relevant example of these wreath products are the double crossed product introduced by Majid  (see \cite{MAJDCP})  which contains, as a particular case, the Drinfel'd double of a Hopf algebra. 

Without leaving the associative setting, the theory of distributive laws was extended to the weak context by Street \cite{RS} and in \cite{BG1} and \cite{BG2} we can find sufficient conditions under which the corresponding weak wreath product, i.e., the  product associated to a weak distributive law, defined as in the first paragraph of this introduction, becomes a weak bimonoid, where the comonoid structure is the tensor product comonoid structure. Also, if the weak bimonoids are weak Hopf monoids, the wreath product is a weak Hopf monoid and, as in the Hopf algebra case, it is possible to describe the Drinfel'd double of a weak Hopf algebra using the wreath product associated to a particular weak distributive law.

Recently, the notion of Hopf quasigroup in a category of vector spaces was introduced  by Klim and Majid in \cite{Majidesfera} in order to understand the structure and relevant properties of the algebraic $7$-sphere, which is the loop of nonzero octonions. Hopf quasigroups are a non-associative generalization of Hopf algebras and they are a particular case  of unital coassociative $H$-bialgebra (see \cite{PI07}) and also of quantum quasigroup (see \cite{SM1} and \cite{SM2}). Hopf quasigroups include as a particular cases the loop algebra for a loop $L$ with the inverse property (see \cite{Majidesfera}, \cite{Bruck}) and also  the enveloping algebra $U(M)$ of a Malcev algebra living in a category of modules over a ring $R$ satisfying suitable conditions  (see \cite{Majidesfera} and \cite{PIS}).  In \cite{our4} the authors proved that  the category of  loops with the inverse property, I.P. loops for short, is equivalent to  the one of  pointed cosemisimple Hopf quasigroups over a given field ${\Bbb F}$ and, if  ${\Bbb F}$ is algebraically closed, the categories of  I.P. loops  and  of cocommutative cosemisimple Hopf quasigroups are equivalent. Also, in \cite{our1} and \cite{our3} was proved that relevant properties of Hopf algebras can be obtained for Hopf quasigroups. Finally, and in a simplified manner, we can say that  Hopf quasigroups unify I.P. loops and Malcev algebras in the same way that Hopf algebras unify groups and Lie algebras. 

In the literature about Hopf quasigroups we can fin several papers devoted to  study  smash products of Hopf quasigroups (see \cite{BrzezJiao1}, \cite{BrzezJiao2}),  twisted smash products of Hopf quasigroups (see \cite{FT1}, \cite{FT}), Hopf quasigroups obtained by the twist double method (see  \cite{FT}), Hopf quasigroups associated to skew pairings  and double cross products of Hopf quasigroups (see \cite{our2}). It is a remarkable fact that in all these products the magma structure is determined by  a morphism $\Psi$  satisfying some conditions that are close to the ones involved in the definition of distributive law. Taking this into consideration,  the main motivation of this paper is to introduce a broad and natural notion of distributive law for Hopf quasigroups that permits to understand the products cited in the first lines of this paragraph with a general point of view.  In order to do it, firstly, in the second section of this paper, for two Hopf quasigropups $A$ and $H$, we introduce the notion of $a$-comonoidal distributive law of $H$ over $A$ and we prove that all the examples listed previously are determined by an $a$-comonoidal distributive law. Secondly, in the last section, we  prove that $A\ot H$ with the corresponding  wreath product, i.e., the  product associated to an $a$-comonoidal distributive law, becomes a  Hopf quasigroup, where the comonoid structure is the one of the tensor product comonoid. As a consequence, we obtain that the double cross product  of two Hopf quasigroups  (see \cite{our2}) is an example of wreath product Hopf quasigroup.  Also is the cross product  associated to a skew pairing (see \cite{our2}).  Moreover,  the Hopf quasigroups defined by the twisted double method in \cite{FT} are examples of wreath product Hopf quasigroups. Finally,  the smash products of Hopf quasigroups in the sense of \cite{BrzezJiao1} are examples of wreath product Hopf quasigroups as well as the twisted smash products defined in \cite{FT1}.

\section{Distributive laws between Hopf quasigroups}

From now on $\mathcal C$ denotes a strict symmetric monoidal category with tensor product $\ot$, unit object $K$ and natural isomorphism of symmetry $c$. Recall that  a monoidal category is a category ${\mathcal C}$ equipped with a tensor product functor $\ot :{\mathcal C}\times {\mathcal C}\rightarrow {\mathcal C}$,  a unit object  $K$ of ${\mathcal C}$ and  a family of natural isomorphisms 
$$a_{M,N,P}:(M\ot N)\ot P\rightarrow M\ot (N\ot P),$$
$$r_{M}:M\ot K\rightarrow M, \;\;\; l_{M}:K\ot M\rightarrow M,$$
in ${\mathcal C}$ (called  associativity, right unit and left unit constraints, respectively) satisfying the Pentagon Axiom and the Triangle Axiom, i.e.,
$$a_{M,N, P\ot Q}\co a_{M\ot N,P,Q}= (id_{M}\ot a_{N,P,Q})\co a_{M,N\ot P,Q}\co (a_{M,N,P}\ot id_{Q}),$$
$$(id_{M}\ot l_{N})\co a_{M,K,N}=r_{M}\ot id_{N},$$
where $id_{X}$ denotes the identity morphism for each object $X$ in ${\mathcal 
C}$ (see \cite{Mac}). A monoidal category is called strict if the associativity, right unit and left unit constraints are identities.  On the other hand, a strict monoidal category ${\mathcal C}$ is symmetric  if it has  a natural family of isomorphisms $c_{M,N}:M\ot N\rightarrow N\ot M$ such that the equalities
$$
c_{M,N\ot P}= (id_{N}\ot c_{M,P})\co (c_{M,N}\ot id_{P}),\;\;
c_{M\ot N, P}= (c_{M,P}\ot id_{N})\co (id_{M}\ot c_{N,P}),\;\; 
$$
$$c_{N,M}\co c_{M,N}=id_{M\ot N},$$
hold for all $M$, $N$ in ${\mathcal C}$.

Taking into account that every non-strict monoidal category is monoidal equivalent to a strict one (see \cite{Christian}), we can assume without loss of generality that the category is strict  and, as a consequence, the results contained in this paper remain valid for every non-strict symmetric monoidal category, what would include for example the categories of vector spaces over a field ${\Bbb F}$, or the one of left modules over a commutative ring $R$. In what follows, for simplicity of notation, given objects $M$, $N$, $P$ in ${\mathcal
C}$ and a morphism $f:M\rightarrow N$, we write $P\ot f$ for
$id_{P}\ot f$ and $f \ot P$ for $f\ot id_{P}$.

A magma in ${\mathcal C}$ is a pair $A=(A,  
\mu_{A})$, where  $A$ is an object in ${\mathcal C}$ and $\mu_{A}:A\otimes A\rightarrow A$ (product) is a morphism in ${\mathcal C}$.
A unital magma  in ${\mathcal C}$ is a triple $A=(A, \eta_A, 
\mu_{A})$, where $(A,  
\mu_{A})$ is a magma in ${\mathcal C}$ and
 $\eta_{A}: K\rightarrow A$ (unit) is a morphism in ${\mathcal C}$ such that 
$\mu_{A}\circ (A\otimes \eta_{A})=id_{A}=\mu_{A}\circ
(\eta_{A}\otimes A)$. A monoid in ${\mathcal C}$ is a unital magma $A=(A, \eta_A, \mu_{A})$ in ${\mathcal C}$ satisfying   $\mu_{A}\circ (A\otimes
\mu_{A})=\mu_{A}\circ (\mu_{A}\otimes A)$, i.e., the product $\mu_{A}$ is associative. Given two unital magmas (monoids) $A$ and $B$,
$f:A\rightarrow B$ is a morphism of unital magmas (monoids) if  $f\circ \eta_{A}= \eta_{B}$ and $\mu_{B}\circ (f\otimes f)=f\circ \mu_{A}$. 

Also,
if $A$, $B$ are unital magmas (monoids) in ${\mathcal C}$, the object $A\otimes
B$ is a unital magma (monoid) in
 ${\mathcal C}$, where $\eta_{A\otimes B}=\eta_{A}\otimes \eta_{B}$
and $\mu_{A\otimes B}=(\mu_{A}\otimes \mu_{B})\circ (A\otimes c_{B,A}\otimes B)$. If $A=(A, \eta_A,  \mu_{A})$ is a unital magma so is $A^{op}=(A, \eta_A,  \mu_{A}\co c_{A,A})$.

A comagma in ${\mathcal C}$ is a pair ${D} = (D, \delta_{D})$, where $D$ is an object in ${\mathcal C}$ and $\delta_{D}:D\rightarrow D\otimes D$ (coproduct) is a morphism in
${\mathcal C}$. A counital comagma in ${\mathcal C}$ is a triple ${D} = (D,
\varepsilon_{D}, \delta_{D})$, where $(D, \delta_{D})$ is a comagma in ${\mathcal
C}$ and $\varepsilon_{D}: D\rightarrow K$ (counit) is a  morphism in
${\mathcal C}$ such that $(\varepsilon_{D}\otimes D)\circ
\delta_{D}= id_{D}=(D\otimes \varepsilon_{D})\circ \delta_{D}$.  A comonoid  in ${\mathcal C}$ is a counital comagma in ${\mathcal C}$ satisfying  $(\delta_{D}\otimes D)\circ \delta_{D}= (D\otimes \delta_{D})\circ \delta_{D}$, i.e., the coproduct $\delta_{D}$ is coassociative. If $D$ and
$E$ are counital comagmas (comonoids) in  ${\mathcal C}$,
$f:D\rightarrow E$ is a  morphism of counital comagmas (comonoids) if $\varepsilon_{E}\circ f
=\varepsilon_{D}$,  and $(f\otimes f)\circ
\delta_{D} =\delta_{E}\circ f$.  

Moreover, if $D$, $E$ are counital comagmas (comonoids) in ${\mathcal C}$,
the object $D\otimes E$ is a counital comagma (comonoid) in ${\mathcal C}$, where
$\varepsilon_{D\otimes E}=\varepsilon_{D}\otimes \varepsilon_{E}$
and $\delta_{D\otimes E}=(D\otimes c_{D,E}\otimes E)\circ(
\delta_{D}\otimes  \delta_{E})$. If $D=(D, \varepsilon_{D}, \delta_{D})$ is a counital comagma so is $D^{cop}=(D, \varepsilon_{D}, c_{D,D}\co\delta_{D})$.

Let $f:B\rightarrow A$ and $g:B\rightarrow A$ be morphisms between a comagma $B$ and a magma $A$. We define the convolution product by $f\ast g=\mu_{A}\circ
(f\otimes g)\circ \delta_{B}$. If $A$ is unital and $B$ counital, we will say that $f$ is convolution invertible if there exists $f^{-1}:B\to A$ such that $f \ast f^{-1}=f^{-1}\ast f=\varepsilon_B\ot \eta_A$. Note that, if $B=K$, we have that $f\ast g=\mu_{A}\circ (f\otimes g)$ and $f$ is convolution invertible if there exists $f^{-1}:K\to A$ such that $f \ast f^{-1}=f^{-1}\ast f= \eta_A$.

\begin{definition}
\label{nbimod}
{\rm A non-associative bimonoid in the category $\mathcal C$ is a unital magma $(H,\eta_H,\mu_H)$ and a comonoid $(H,\varepsilon_H,\delta_H)$ such that $\varepsilon_H$ and $\delta_H$ are morphisms of unital magmas (equivalently, $\eta_H$ and $\mu_H$ are morphisms of counital comagmas). Then the following identities hold:
\begin{equation}
\label{eta-eps}
\varepsilon_{H}\co \eta_{H}=id_{K},
\end{equation}
\begin{equation}
\label{mu-eps}
\varepsilon_{H}\co \mu_{H}=\varepsilon_{H}\ot \varepsilon_{H},
\end{equation}
\begin{equation}
\label{delta-eta}
\delta_{H}\co \eta_{H}=\eta_{H}\ot \eta_{H},
\end{equation}
\begin{equation}
\label{delta-mu}
\delta_{H}\co \mu_{H}=(\mu_{H}\ot \mu_{H})\co \delta_{H\ot H}.
\end{equation}

}
\end{definition}

A non-associative bimonoid is called cocommutative if $\delta_{H}=c_{H,H}\co \delta_{H}$, i.e., $H=H^{cop}$ as comonoids.

Now we introduce the notion of  Hopf quasigroup.

\begin{definition}
\label{Hopfqg} {\rm A Hopf quasigroup $H$ in
${\mathcal C}$ is a non-associative bimonoid such that there exists a morphism $\lambda_{H}:H\rightarrow H$ in ${\mathcal C}$ (called the  antipode of $H$) satisfying
\begin{equation}
\label{lH}
\mu_H\circ (\lambda_H\ot \mu_H)\circ (\delta_H\ot H)=
\varepsilon_H\ot H=
\mu_H\circ (H\ot \mu_H)\circ (H\ot \lambda_H\ot H)\circ (\delta_H\ot H)
\end{equation}
and
\begin{equation}
\label{rH}
\mu_H\circ (\mu_H\ot H)\circ (H\ot \lambda_H\ot H)\circ (H\ot \delta_H)=
H\ot \varepsilon_H=
\mu_H\circ(\mu_H\ot \lambda_H)\circ (H\ot \delta_H). 
\end{equation}

Note that composing with $H\ot \eta_{H}$ in (\ref{lH}) we obtain that 
\begin{equation}
\label{1-conv}
\lambda_{H}\ast id_{H}=\varepsilon_H\ot \eta_H,  
\end{equation}
and composing with $ \eta_{H}\ot H$ in (\ref{rH}) we obtain 
\begin{equation}
\label{2-conv}
id_{H}\ast \lambda_H=\varepsilon_H\ot \eta_H. 
\end{equation}

Therefore, $\lambda_{H}$ is convolution invertible and $\lambda_{H}^{-1}=id_{H}$.

}
\end{definition}

The above definition is the monoidal version of the notion of Hopf quasigroup (also called non-associative Hopf algebra with the inverse property, or non-associative IP Hopf algebra) introduced in \cite{Majidesfera} (in this case ${\mathcal C}$ is the category of vector spaces over a field ${\Bbb F}$).   Note that a  Hopf quasigroup is associative, i.e., $H$ is a monoid,  if and only if  it is a Hopf algebra.

If $H$ is a Hopf quasigroup in ${\mathcal C}$  we know that the antipode
$\lambda_{H}$ is unique because, if  $s_{H}:H\rightarrow H$ is a new morphism satisfying (\ref{lH}) and (\ref{rH}), then
$$s_{H}= (s_{H}\ast id_{H})\ast s_{H} = (\lambda_{H}\ast id_{H})\ast s_{H} = \mu_{H}\co (\mu_{H}\ot s_{H})\co (\lambda_{H}\ot \delta_{H})\co \delta_{H}= \lambda_{H}.$$ 

Also, by \cite[Proposition 4.2]{Majidesfera}, the antipode is antimultiplicative and anticomultiplicative, i.e.,
\begin{equation} 
\label{antimu}
\lambda_{H}\co\mu_{H}=\mu_{H}\circ
(\lambda_{H}\ot\lambda_{H})\circ
c_{H,H},
\end{equation}
\begin{equation} 
\label{anticm}
\delta_{H}\co \lambda_{H}=
(\lambda_{H}\ot\lambda_{H})\circ
c_{H,H}\co \delta_{H}
\end{equation}
and leaves the unit and the counit invariable (see \cite{LP}):
\begin{equation} 
\label{inv}
\lambda_{H}\co\eta_{H}=\eta_{H},\;\;\;\varepsilon_{H}\co \lambda_{H}=\varepsilon_{H}. 
\end{equation}

A morphism between Hopf quasigroups $H$ and $B$ is a morphism $f:H\rightarrow B$ of unital magmas and comonoids. Then (see Lemma 1.4 of \cite{our1}) the equality
\begin{equation}
\label{antipode-morphism}
\lambda_B\co f=f\co \lambda_H
\end{equation}
holds.

\begin{example}
\label{exloop}
{\rm  A loop $(L, \cdot,\diagup, \diagdown)$ is a quadruple where $L$ is a set, $\cdot$ (multiplication),  $\diagup$ (right division) and
$\diagdown$ (left division ) are binary operations, satisfying the identities
$$
 v\diagdown(v \cdot u) = u,\;\;
 u = (u \cdot v)\diagup v,\;\;
v\cdot  (v\diagdown u) = u,\;\;
 u = (u\diagup v) \cdot  v,
$$
and such that it contains an  identity element $e_L$ (i.e., $e_L\cdot x=x=x\cdot e_L$ hold for all $x$ in $L$). In what follows multiplications on $L$ will be expressed by juxtaposition. 

Let $R$ be a commutative ring and $L$ a loop. Then, the loop algebra 
$$RL=\bigoplus_{u\in L}Ru$$
is a cocommutative non-associative bimonoid with unit $\eta_{RL}(1_{R})=e_{L}$,  product  defined by
linear extension of the one defined in $L$, $\delta_{RL}(u)=u\ot u $ and $\varepsilon_{RL}(u)=1_{R}$ 
 on the basis elements (see \cite{PI07}).
 
Let $L$ be a loop. If 
for every element $u\in L$ there exists an element $u^{-1}\in
L$ (the inverse of $u$) such that the equalities $ u^{-1}(uv)=v=(vu)u^{-1}$
hold for  every $v\in L$, we will say that $L$ is a loop with the
inverse property  (for brevity an IP loop). As a consequence, it is easy to show that , if $L$ is an IP loop, 
for all $u\in L$ the element $u^{-1}$ is unique and $ u^{-1}u=e_{L}=uu^{-1}$.  Moreover, $ (uv)^{-1}=v^{-1}u^{-1} $ holds for any pair of elements $u,v\in L$.

Now let $R$ be a commutative ring and let $L$ be and IP loop. Then, by
\cite[Proposition 4.7]{Majidesfera}, the non-associative bimonoid  $RL$ is a cocommutative Hopf quasigroup 
where the antipode is defined by the linear extension of the map $\lambda_{RL}(u)=u^{-1}$.

}
\end{example}

\begin{example}
\label{exsabinin}
{\rm
Let $R$ be a commutative ring  with $\frac{1}{2}$ and $\frac{1}{3}$ in $R$. A Malcev algebra  $(M,[,])$ over $R$ is a free $R$-module $M$ equipped with a bilinear and anticommutative
operation [ , ] such that
\[[J(a, b, c), a] = J(a, b, [a, c]),\]
where $J(a, b, c) = [[a, b], c] - [[a, c], b] - [a, [b, c]]$ denotes the Jacobian in $a$, $b$, $c$ (see \cite{PIS}). Then,  every Lie algebra is a Malcev algebra with $J=0$.  The universal enveloping algebra $U(M)$ can be provided with a Hopf quasigroup structure \cite[Proposition 4.8]{Majidesfera}.

}
\end{example}

\begin{definition}
\label{dl} {\rm Le $H$, $A$ be Hopf quasigroups. A morphism $\Psi:H\ot A\rightarrow A\ot H$ is said to be a distributive law of $H$ over $A$ if the following identities
\begin{equation}
\label{dl1}
\Psi\co (H\ot \mu_{A})\co (\lambda_{H}\ot \lambda_{A}\ot A)=(\mu_{A}\ot H)\co (A\ot \Psi)\co (\Psi\ot A)\co (\lambda_{H}\ot \lambda_{A}\ot A),
\end{equation}
\begin{equation}
\label{dl2}
\Psi\co (\mu_{H}\ot A)\co (H\ot \lambda_{H}\ot \lambda_{A})=(A\ot \mu_{H})\co ( \Psi\ot H)\co (H\ot \Psi)\co (H\ot \lambda_{H}\ot \lambda_{A}),
\end{equation}
\begin{equation}
\label{dl3}
\Psi\co (H\ot \eta_{A})=\eta_{A}\ot H,
\end{equation}
\begin{equation}
\label{dl4}
\Psi\co (\eta_{H}\ot A)=A\ot \eta_{H}, 
\end{equation}
hold.
}
\end{definition}

If   the antipodes of $H$ and $A$ are isomorphisms, the identities (\ref{dl1}) and (\ref{dl2}) are equivalent to 
\begin{equation}
\label{dl1-1}
\Psi\co (H\ot \mu_{A})=(\mu_{A}\ot H)\co (A\ot \Psi)\co (\Psi\ot A),
\end{equation}
\begin{equation}
\label{dl2-1}
\Psi\co (\mu_{H}\ot A)=(A\ot \mu_{H})\co ( \Psi\ot H)\co (H\ot \Psi),
\end{equation}
respectively. Then, in this case, the conditions of the definition of distributive law for Hopf quasigroups are the ones that we can find in the classical definition of distributive law between monoids, i.e., $\Psi$ is compatible with the unit and the product of $A$ and $H$. 

\begin{definition}
\label{cdl} {\rm Le $H$, $A$ be Hopf quasigroups and let  $\Psi:H\ot A\rightarrow A\ot H$  be a distributive law of $H$ over $A$. The distributive law $\Psi$ is said to be comonoidal if it is a comonoid morphism, i.e., the following identities
\begin{equation}
\label{cdl1}
\delta_{A\ot H}\co \Psi=(\Psi \ot \Psi)\co \delta_{H\ot A},
\end{equation}
\begin{equation}
\label{cdl2}
(\varepsilon_{A}\ot \varepsilon_{H})\co\Psi=\varepsilon_{H}\ot \varepsilon_{A}, 
\end{equation}
hold.
}
\end{definition}

\begin{definition}
\label{cdl} {\rm Le $H$, $A$ be Hopf quasigroups and let  $\Psi:H\ot A\rightarrow A\ot H$  be a comonoidal distributive law of $H$ over $A$. We will say that $\Psi$ is  an $a$-comonoidal distributive law of $H$ over $A$  if the following identities
\begin{equation}
\label{adl1}
(A\ot \mu_{H})\co (\Psi\ot \mu_{H})\co (H\ot \Psi\ot H)\co (((\lambda_{H}\ot H)\co \delta_{H})\ot A\ot H)=\varepsilon_{H}\ot A\ot H,
\end{equation}
\begin{equation}
\label{adl2}
(A\ot \mu_{H})\co (\Psi\ot \mu_{H})\co (H\ot \Psi\ot H)\co (((H\ot \lambda_{H})\co \delta_{H})\ot A\ot H)=\varepsilon_{H}\ot A\ot H,
\end{equation}
\begin{equation}
\label{adl3}
(\mu_{A}\ot H)\co (\mu_{A}\ot \Psi)\co (A\ot \Psi\ot A)\co (A\ot H\ot (( \lambda_{A}\ot A)\co \delta_{A}))=A\ot H\ot \varepsilon_{A}, 
\end{equation}
\begin{equation}
\label{adl4}
(\mu_{A}\ot H)\co (\mu_{A}\ot \Psi)\co (A\ot \Psi\ot A)\co (A\ot H\ot (( A\ot \lambda_{A})\co \delta_{A}))=A\ot H\ot \varepsilon_{A},
\end{equation} 
hold.
}
\end{definition}

Note that, if $H$ and $A$ are Hopf algebras and $\Psi:H\ot A\rightarrow A\ot H$ is a distributive law between the monoids $H$ and $A$, the equalities (\ref{adl1}), (\ref{adl2}), (\ref{adl3}) and (\ref{adl4}) always hold.

\begin{example}
\label{1}
{\rm In this example we will show that the theory of double cross products of  Hopf quasigroups, introduced in \cite{our2}, produces examples of $a$-comonoidal distributive laws. 

Let  $H$ be a  Hopf quasigroup. The pair $(M,\varphi_{M})$ is said to be a  left $H$-quasimodule if $M$ is an object in ${\mathcal C}$ and $\varphi_{M}:H\ot M\rightarrow M$ is a morphism in ${\mathcal C}$ (called the action) satisfying
\begin{equation}
\label{uq}
\varphi_{M}\circ(\eta_{H}\ot M)=id_{M}
\end{equation}
and
\begin{equation}
\label{pq}
\varphi_{M}\co (H\ot\varphi_{M})\co  (((H\ot \lambda_{H})\co \delta_H)\ot M)=\varepsilon_H\ot M=\varphi_{M}\co (\lambda_H\ot \varphi_M)\co (\delta_H\ot M).
\end{equation}

Given two left ${H}$-quasimodules $(M,\varphi_{M})$, $(N,\varphi_{N})$ and a morphism  $f:M\rightarrow N$ in ${\mathcal C}$, we will say that $f$ is a morphism of left
$H$-quasimodules if 
\begin{equation}
\label{mor}
\varphi_{N}\co(H\ot f)=f\co\varphi_{M}.
\end{equation}

We denote the category of left $H$-quasimodules by $_{H}{\mathcal {QC}}$. It is easy to prove that, if $(M,\varphi_{M})$ and $(N,\varphi_{N})$ are left $H$-quasimodules, the tensor product 
 $M\ot N$ is a left $H$-quasimodule with the diagonal action 
 $$\varphi_{M\otimes
N}=(\varphi_{M}\ot\varphi_{N})\circ
(H\ot c_{H,M}\ot N)\co(\delta_{H}\ot M\ot N).$$

This makes the category of left $H$-quasimodules into a strict monoidal category $(_{H}{\mathcal {QC}}, \ot , K)$ (see \cite[Remark 3.3]{BrzezJiao1}).

We will say that a unital magma $A$ is a left $H$-quasimodule magma  if it is a left $H$-quasimodule with action $\varphi_{A}: H\ot A\rightarrow A$ and the following equalities 
\begin{equation}
\label{etaq}
\varphi_{A}\circ(H\ot \eta_A)=\varepsilon_H\ot \eta_A,
\end{equation}
\begin{equation}
\label{muq}
\mu_A\co \varphi_{A\ot A}=\varphi_A\co (H\ot \mu_A),
\end{equation}
hold, i.e., $\eta_{A}$ and $\mu_{A}$ are quasimodule morphisms.

A comonoid $A$ is a left $H$-quasimodule comonoid if it is a left $H$-quasimodule with action $\varphi_{A}$ and 
\begin{equation}
\label{eq}
\varepsilon_A\co \varphi_{A}=\varepsilon_H\ot \varepsilon_A,
\end{equation}
\begin{equation}
\label{dq}
\delta_A\co \varphi_{A}=\varphi_{A\ot A}\co (H\ot\delta_A),
\end{equation}
hold, i.e., $\varepsilon_{A}$ and $\delta_{A}$ are quasimodule morphisms.

Replacing (\ref{pq}) by the equality 
\begin{equation}
\label{pqmod}
\varphi_{M}\co(H\ot\varphi_{M})=\varphi_M\co (\mu_H\ot M),
\end{equation}
we have the definition of left $H$-module  and the ones of left $H$-module magma and comonoid. Note that the pair $(H, \mu_H)$ is not an $H$-module but it is an  $H$-quasimodule. Morphisms between left $H$-modules are defined as for $H$-quasimodules and we denote the category of left $H$-modules by $\;_{H}{\mathcal C}$. Obviously we have similar definitions for the right side.

By \cite[Proposition 5.3]{our2} and  its right version,  in \cite[Corollary 5.6]{our2} we prove that, if $A$, $H$ are  Hopf quasigroups,  $(A, \varphi_A)$ is a left $H$-module comonoid, $(H, \phi_H)$  is a right $A$-module comonoid and 
$$\Psi=(\varphi_A\ot \phi_H)\co \delta_{H\ot A},$$ 
the following assertions are equivalent:

\begin{itemize}
\item [(i)] The double cross product $A\bowtie H$ built on the object $A\ot H$ with product
$$\mu_{A\bowtie H}=(\mu_A\ot \mu_H)\co (A\ot \Psi\ot H)$$
and tensor product unit, counit and coproduct, is a Hopf quasigroup with antipode
$$\lambda_{A\bowtie H}=\Psi\co (\lambda_H\ot \lambda_A)\co c_{A,H}.$$

\item [(ii)] The  equalities

 \begin{equation}\label{d1}
\varphi_A\co (H\ot \eta_A)=\varepsilon_H\ot \eta_A,
 \end{equation}
\begin{equation}\label{d2}
\phi_H\co (\eta_H\ot A)=\eta_H\ot \varepsilon_A,
 \end{equation}
\begin{equation}\label{d3}
(\phi_H\ot \varphi_A)\co \delta_{H\ot A}=c_{A,H}\co \Psi,
 \end{equation}
\begin{equation}\label{d4}
\varphi_A\co (H\ot \mu_A)\co (\lambda_H\ot \lambda_A\ot A)=\mu_A\co(A\ot \varphi_A)\co ((\Psi\co (\lambda_H\ot \lambda_A))\ot A),
 \end{equation}
\begin{equation}\label{d5}
\mu_H\co(\phi_H\ot \mu_H)\co (\lambda_H\ot \Psi \ot H)\co  (\delta_H\ot A\ot H)=\varepsilon_H\ot \varepsilon_A\ot H,
 \end{equation}
\begin{equation}\label{d6}
\mu_H\co(\phi_H\ot \mu_H)\co (H\ot \Psi\ot H)\co  (((H\ot \lambda_H)\co\delta_H)\ot A\ot H)=\varepsilon_H\ot \varepsilon_A\ot H,
 \end{equation}
 \begin{equation}\label{d7}
\phi_H\co (\mu_{H}\ot A)\co (H\ot\lambda_H\ot \lambda_A)=\mu_H\co(\phi_{H}\ot H)\co (H\ot (\Psi\co (\lambda_H\ot \lambda_A))),
 \end{equation}
\begin{equation}\label{d8}
\mu_A\co(\mu_{A}\ot \varphi_A)\co (A\ot \Psi\ot \lambda_{A})\co (A\ot H\ot \delta_{A})=A\ot \varepsilon_H\ot \varepsilon_A,
 \end{equation}
\begin{equation}\label{d9}
\mu_{A}\co (\mu_{A}\ot \varphi_{A})\co (A\ot \Psi\ot A)\co (A\ot H\ot ((\lambda_{A}\ot A)\co \delta_{A}))=A\ot \varepsilon_H\ot \varepsilon_A,
 \end{equation}
hold.

\end{itemize}

If the previous equalities hold, we can prove that $\Psi$ is an example of $a$-comonoidal distributive law of $H$ over $A$. Indeed:  First,   
\begin{itemize}
\item[ ]$\hspace{0.38cm}(\mu_{A}\ot H)\co (A\ot \Psi)\co (\Psi\ot A)\co (\lambda_{H}\ot \lambda_{A}\ot A)  $

\item[ ]$=(\mu_{A}\ot H)\co (A\ot ((\varphi_{A}\ot \phi_{H})\co (H\ot c_{H,A}\ot A)))\co (\varphi_{A}\ot (\delta_{H}\co \phi_{H})\ot A\ot A)\co ((\delta_{A\ot H}$
\item[ ]$\hspace{0.38cm}\co (\lambda_{A}\ot \lambda_{H}))\ot \delta_{A})$ {\scriptsize ({\blue by definition})}

\item[ ]$=(\mu_{A}\ot H)\co (A\ot ((\varphi_{A}\ot \phi_{H})\co (H\ot c_{H,A}\ot A)))\co (\varphi_{A}\ot ((\phi_{H}\ot \phi_{H})\co \delta_{H\ot A})\ot A\ot A)$
\item[ ]$\hspace{0.38cm}\co ((\delta_{A\ot H}\co (\lambda_{A}\ot \lambda_{H}))\ot \delta_{A})$ {\scriptsize ({\blue by the condition of of right $A$-module comonoid for $(H,\phi_{H})$})}

\item[ ]$=(\mu_{A}\ot H)\co (A\ot ((\varphi_{A}\ot \phi_{H})\co (H\ot c_{H,A}\ot A)))\co (\Psi\ot H\ot A\ot \mu_{A})\co (H\ot A\ot H\ot c_{A,A}\ot A)$
\item[ ]$\hspace{0.38cm}\co ((\delta_{A\ot H}\co (\lambda_{A}\ot \lambda_{H}))\ot \delta_{A})$  {\scriptsize ({\blue by the condition of of right $A$-module for $(H,\phi_{H})$, naturality of $c$, and }}
\item[ ]$\hspace{0.38cm}$   {\scriptsize {\blue   coassociativity of $\delta_{A}$ and $\delta_{H}$})}

\item[ ]$= ((\mu_A\co(A\ot \varphi_A)\co ((\Psi\co (\lambda_H\ot \lambda_A))\ot A))\ot \phi_{H})\co (H\ot A\ot c_{H,A}\ot \mu_{A})\co (H\ot c_{H,A}\ot c_{A,A}\ot A)$
\item[ ]$\hspace{0.38cm}\co ((c_{H,H}\co (\lambda_{H}\ot H)\co \delta_{H})\ot (c_{A,A}\co (\lambda_{A}\ot A)\co \delta_{A})\ot \delta_{A})$  {\scriptsize ({\blue by (\ref{anticm}) for $\lambda_{A}$ and $\lambda_{H}$ and naturality}}
\item[ ]$\hspace{0.38cm}$ {\scriptsize {\blue of $c$})}

\item[ ]$= ((\varphi_A\co (H\ot \mu_A)\co (\lambda_H\ot \lambda_A\ot A))\ot A))\ot \phi_{H})\co (H\ot A\ot c_{H,A}\ot \mu_{A})\co (H\ot c_{H,A}\ot c_{A,A}\ot A)$
\item[ ]$\hspace{0.38cm}\co ((c_{H,H}\co (\lambda_{H}\ot H)\co \delta_{H})\ot (c_{A,A}\co (\lambda_{A}\ot A)\co \delta_{A})\ot \delta_{A})$  {\scriptsize ({\blue by(\ref{d4})})}

\item[ ]$=(\varphi_{A}\ot H)\co (H\ot \mu_{A}\ot \phi_{H})\ot (H\ot A\ot c_{H,A}\ot \mu_{A})\co (H\ot A\ot H\ot c_{A,A}\ot A)\co ((\delta_{A\ot H}$ 
\item[ ]$\hspace{0.38cm}\co (\lambda_{A}\ot \lambda_{H}))\ot \delta_{A})$ {\scriptsize  ({\blue  by (\ref{anticm}) for $\lambda_{A}$ and $\lambda_{H}$ and naturality of $c$})}

\item[ ]$=(\varphi_{A}\ot \phi_{H})\co (H\ot c_{H,A}\ot A)\co (\delta_{H}\ot ((\mu_{A}\ot \mu_{A})\co \delta_{A\ot A}))\co (\lambda_{H}\ot \lambda_{A}\ot A)$  {\scriptsize  ({\blue by naturality}}
\item[ ]$\hspace{0.38cm}$ {\scriptsize  {\blue  of $c$})}

\item[ ]$= \Psi\co (H\ot \mu_{A})\co (\lambda_{H}\ot \lambda_{A}\ot A)$ {\scriptsize  ({\blue  by (\ref{delta-mu}) for $A$})}
\end{itemize}
and then (\ref{dl1}) holds. Second,  (\ref{dl2}) holds because 

\begin{itemize}
\item[ ]$\hspace{0.38cm}(A\ot \mu_{H})\co ( \Psi\ot H)\co (H\ot \Psi)\co (H\ot \lambda_{H}\ot \lambda_{A})  $

\item[ ]$=(A\ot \mu_{H})\co ( ((\varphi_{A}\ot \phi_{H})\co (H\ot c_{H,A}\ot A))\ot H)\co (H\ot H\ot (\delta_{A}\co \varphi_{A})\ot \phi_{H})\co (\delta_{H}\ot (\delta_{H\ot A}$
\item[ ]$\hspace{0.38cm}\co (\lambda_{H}\ot \lambda_{A})))$ {\scriptsize ({\blue by definition})}

\item[ ]$=(A\ot \mu_{H})\co ( ((\varphi_{A}\ot \phi_{H})\co (H\ot c_{H,A}\ot A))\ot H)\co (H\ot H\ot ((\varphi_{A}\ot \varphi_{A})\co \delta_{H\ot A})\ot \phi_{H})$
\item[ ]$\hspace{0.38cm}\co (\delta_{H}\ot (\delta_{H\ot A}\co (\lambda_{H}\ot \lambda_{A})))$ {\scriptsize ({\blue by (\ref{dq})})}

\item[ ]$=(A\ot \mu_{H})\co  ( ((\varphi_{A}\ot \phi_{H})\co (H\ot c_{H,A}\ot A))\ot H)\co (\mu_{H}\ot c_{H,A}\ot \Psi)\co (H\ot c_{H,H}\ot A\ot H\ot A)$
\item[ ]$\hspace{0.38cm}\co (\delta_{H}\ot (\delta_{H\ot A}\co (\lambda_{H}\ot \lambda_{A})))$  {\scriptsize ({\blue by the condition of of left $H$-module for $(A,\varphi_{A})$, naturality of $c$, and }}
\item[ ]$\hspace{0.38cm}$   {\scriptsize {\blue  coassociativity of $\delta_{A}$ and $\delta_{H}$})}

\item[ ]$= (\varphi_{A}\ot (\mu_H\co(\phi_{H}\ot H)\co (H\ot (\Psi\co (\lambda_H\ot \lambda_A))))\co (\mu_{H}\ot c_{H,A}\ot H\ot A)\co (H\ot c_{H,H}\ot c_{H,A}\ot A)$
\item[ ]$\hspace{0.38cm}\co (\delta_{H}\ot (c_{H,H}\co (H\ot \lambda_{H})\co \delta_{H})\ot (c_{A,A}\co (A\ot \lambda_{A})\co \delta_{A}))$  {\scriptsize ({\blue by (\ref{anticm}) for $\lambda_{A}$ and $\lambda_{H}$ and   naturality}}
\item[ ]$\hspace{0.38cm}$ {\scriptsize {\blue of $c$})}

\item[ ]$= (\varphi_{A}\ot (\phi_H\co (\mu_{H}\ot A)\co (H\ot\lambda_H\ot \lambda_A)))\co (\mu_{H}\ot c_{H,A}\ot H\ot A)\co (H\ot c_{H,H}\ot c_{H,A}\ot A)$
\item[ ]$\hspace{0.38cm}\co (\delta_{H}\ot (c_{H,H}\co (H\ot \lambda_{H})\co \delta_{H})\ot (c_{A,A}\co (A\ot \lambda_{A})\co \delta_{A}))$  {\scriptsize ({\blue by(\ref{d7})})}

\item[ ]$=(A\ot \phi_{H})\co (\varphi_{A}\ot \mu_{H}\ot A)\ot (\mu_{H}\ot c_{H,A}\ot H\ot A)\co (H\ot c_{H,H}\ot A\ot H\ot A)\co (\delta_{H}\ot (\delta_{H\ot A}$ 
\item[ ]$\hspace{0.38cm}\co (\lambda_{H}\ot \lambda_{A})))$ {\scriptsize  ({\blue  by (\ref{anticm}) for $\lambda_{A}$ and $\lambda_{H}$ and naturality of $c$})}

\item[ ]$=(\varphi_{A}\ot \phi_{H})\co (H\ot c_{H,A}\ot A)\co (((\mu_{H}\ot \mu_{H})\co \delta_{H\ot H})\ot \delta_{A})\co (H\ot \lambda_{H}\ot \lambda_{A})$ {\scriptsize  ({\blue by naturality}}
\item[ ]$\hspace{0.38cm}$ {\scriptsize  {\blue  of $c$})}

\item[ ]$= \Psi\co (\mu_{H}\ot A)\co (H\ot \lambda_{H}\ot \lambda_{A})$ {\scriptsize  ({\blue  by (\ref{delta-mu}) for $H$}).}
\end{itemize}

On the other hand,  (\ref{dl3}) follows by 
$$\Psi\co (H\ot \eta_{A})=((\varphi_{A}\co (H\ot \eta_{A}))\ot (\phi_{H}\co (H\ot \eta_{A})))\co \delta_{H}=\eta_{A}\ot ((\varepsilon_{H}\ot H)\co \delta_{H})=\eta_{A}\ot H,$$
where the first identity is a consequence of (\ref{delta-eta}) for $A$ and the naturality of $c$, in the second one we used (\ref{d1}) and the right $A$-module condition for $(H, \phi_{H})$ and, finally, the third one follows by the properties of the counit of $H$.  Similarly, by (\ref{delta-eta}) for $H$, the naturality of $c$,  (\ref{d2}), the left $A$-module condition for $(A, \varphi_{H})$ and  the properties of the counit of $A$ we obtain  (\ref{dl4}). 

Therefore $\Psi$ is s a distributive law of $H$ over $A$. Also,  $\Psi$ is comonoidal because, in one hand,  (\ref{cdl1}) follows by 

\begin{itemize}
\item[ ]$\hspace{0.38cm}\delta_{A\ot H}\co \Psi $

\item[ ]$=(A\ot c_{A,H}\ot H)\co ((\delta_{A}\co \varphi_{A})\ot (\delta_{H}\co \phi_{H}))\co \delta_{H\ot A}$  {\scriptsize  ({\blue by definition})}

\item[ ]$=(A\ot c_{A,H}\ot H)\co (\varphi_{A}\ot \varphi_{A}\ot \phi_{H}\ot \phi_{H})\co (\delta_{H\ot A}\ot \delta_{H\ot A})\co \delta_{H\ot A}$  {\scriptsize  ({\blue by (\ref{dq}) and the similar}}
\item[ ]$\hspace{0.38cm}$ {\scriptsize  {\blue  property for $\Phi_{H}$})}

\item[ ]$=(\varphi_{A}\ot (c_{A,H}\co \Psi)\ot \phi_{H})\co (H\ot c_{H,A}\ot c_{H,A}\ot A)\co (\delta_{H}\ot A\ot H\ot \delta_{A})\co \delta_{H\ot A}$  {\scriptsize  ({\blue by naturality}}
\item[ ]$\hspace{0.38cm}$ {\scriptsize  {\blue of $c$ and  coassociativity of $\delta_{A}$ and $\delta_{H}$})}

\item[ ]$=(\varphi_{A}\ot ((\phi_{H}\ot \varphi_{A})\co \delta_{H\ot A})\ot \phi_{H})\co (H\ot c_{H,A}\ot c_{H,A}\ot A)\co (\delta_{H}\ot A\ot H\ot \delta_{A})\co \delta_{H\ot A}$  
\item[ ]$\hspace{0.38cm}$ {\scriptsize  ({\blue by (\ref{d3})})}

\item[ ]$= (\Psi\ot \Psi)\co \delta_{H\ot A}$ {\scriptsize  ({\blue  by naturality of $c$ and  coassociativity of $\delta_{A}$ and $\delta_{H}$}) }
\end{itemize}
and, on the other hand, by (\ref{eq}) for $\varphi_{A}$, the similar identity for $\Phi_{H}$ and the properties of the counits $\varepsilon_{H}$, $\varepsilon_{A}$ we obtain (\ref{cdl2}). 

Finally, we prove that $\Psi$ is $a$-comonoidal. First, to demonstrate  the identities (\ref{adl1}), (\ref{adl2}), (\ref{adl3}) and (\ref{adl3}), we need to show that 
\begin{equation}
\label{hp}
(A\ot \mu_{H})\co (\Psi \ot \mu_{H})\co (H\ot \Psi\ot H) 
\end{equation}
$$=((\varphi_{A}\co (\mu_{H}\ot A))\ot (\mu_{H}\co (\phi_{H}\ot \mu_{H})\co (H\ot \Psi\ot H)))\co (\delta_{H\ot H\ot A}\ot H),$$
\begin{equation}
\label{ap}
(\mu_{A}\ot H)\co (\mu_{A}\ot \Psi)\co (H\ot \Psi\ot A) 
\end{equation}
$$=((\mu_{A}\co (\mu_{A}\ot \varphi_{A})\co (H\ot \Psi\ot A))\ot (\phi_{H}\co (H\ot \mu_{A})))\co (A\ot \delta_{H\ot A\ot A}), $$
hold. Indeed, the equality (\ref{hp}) follows by 

\begin{itemize}
\item[ ]$\hspace{0.38cm} (A\ot \mu_{H})\co (\Psi \ot \mu_{H})\co (H\ot \Psi\ot H) $

\item[ ]$=(\varphi_{A}\ot (\mu_{H}\co (\phi_{H}\ot H)))\co (H\ot c_{H,A}\ot A\ot H)\co (\delta_{H}\ot ((\delta_{A}\co \varphi_{A})\ot (\mu_{H}\co (\phi_{H}\ot H)))$ 
\item[ ]$\hspace{0.38cm}\co  (H\ot \delta_{H\ot A}\ot H)$ {\scriptsize  ({\blue by definition})} 

\item[ ]$=(\varphi_{A}\ot (\mu_{H}\co (\phi_{H}\ot H)))\co (H\ot c_{H,A}\ot A\ot H)\co (\delta_{H}\ot (((\varphi_{A}\ot \varphi_{A})\co \delta_{H\ot A})\ot (\mu_{H}\co (\phi_{H}\ot H))) $  
\item[ ]$\hspace{0.38cm} \co (H\ot \delta_{H\ot A}\ot H)$ {\scriptsize  ({\blue by (\ref{dq})})}

\item[ ]$=(\varphi_{A}\ot (\mu_{H}\co (\phi_{H}\ot H)))\co (H\ot c_{H,A}\ot A\ot \mu_{H})\co (\delta_{H}\ot ((\varphi_{A}\ot \Psi)\co \delta_{H\ot A})\ot H)${\scriptsize  ({\blue by }} 
\item[ ]$\hspace{0.38cm}$   {\scriptsize  {\blue coassociativity of $\delta_{H\ot A}$}}

\item[ ]$= ((\varphi_{A}\co (\mu_{H}\ot A))\ot (\mu_{H}\co (\phi_{H}\ot \mu_{H})\co (H\ot \Psi\ot H)))\co (\delta_{H\ot H\ot A}\ot H)$ {\scriptsize  ({\blue  by naturality of $c$ }}
\item[ ]$\hspace{0.38cm}$   {\scriptsize  {\blue and  the condition of left $H$-module for $(A,\varphi_{A})$}}
\end{itemize}
and similarly, using the condition of right $A$-module comonoid for $(H, \phi_{H})$ we obtain (\ref{ap}). 

The identity  (\ref{adl1}) follows by 

\begin{itemize}
\item[ ]$\hspace{0.38cm} (A\ot \mu_{H})\co (\Psi\ot \mu_{H})\co (H\ot \Psi\ot H)\co (((\lambda_{H}\ot H)\co \delta_{H})\ot A\ot H)$

\item[ ]$=((\varphi_{A}\co (\mu_{H}\ot A))\ot (\mu_{H}\co (\phi_{H}\ot \mu_{H})\co (H\ot \Psi\ot H)))\co (\delta_{H\ot H\ot A}\ot H)\co   (((\lambda_{H}\ot H)$
\item[ ]$\hspace{0.38cm}\co \delta_{H})\ot A\ot H)${\scriptsize ({\blue by (\ref{hp})})} 

\item[ ]$=((\varphi_{A}\co ((\lambda_{H}\ast id_{H})\ot A))\ot (\mu_{H}\co (\phi_{H}\ot \mu_{H})\co (H\ot \Psi\ot H)))\co (H\ot c_{H,A}\ot H\ot A\ot H)$  
\item[ ]$\hspace{0.38cm} \co (c_{H,H}\ot c_{H,A}\ot A\ot H)\co    (((\lambda_{H}\ot \delta_{H})\co \delta_{H})\ot \delta_{A}\ot H)$ {\scriptsize  ({\blue by coassociativity of $\delta_{H}$, naturality }}
\item[ ]$\hspace{0.38cm}${\scriptsize  {\blue of $c$ and (\ref{anticm}) for $\lambda_{H}$})}

\item[ ]$= (A\ot (\mu_{H}\co (\phi_{H}\ot \mu_{H})\co (H\ot \Psi\ot H)))\co ( c_{H,A}\ot H\ot A\ot H) \co  (H\ot  c_{H,A}\ot H\ot A)\co (((\lambda_{H}\ot H)$
\item[ ]$\hspace{0.38cm} \co \delta_{H})\ot \delta_{A}\ot H)$ {\scriptsize  ({\blue by (\ref{1-conv}) for $H$, counit properties, naturality of $c$ and the condition of left $H$-module for $A$})} 

\item[ ]$=(A\ot ( \mu_H\co(\phi_H\ot \mu_H)\co (\lambda_H\ot \Psi \ot H)\co  (\delta_H\ot A\ot H)))\co (c_{A,H}\ot A\ot H)\co (H\ot \delta_{A}\ot H)$ 
\item[ ]$\hspace{0.38cm}$ {\scriptsize  ({\blue  by naturality of $c$})}

\item[ ]$=(A\ot \varepsilon_{H}\ot \varepsilon_{A}\ot H)\co (c_{A,H}\ot A\ot H)\co (H\ot \delta_{A}\ot H)$ {\scriptsize  ({\blue by (\ref{d5})})}

\item[ ]$=\varepsilon_{H}\ot A\ot H$  {\scriptsize  ({\blue by naturality of $c$ and counit properties}).}

\end{itemize}

Also, 

\begin{itemize}
\item[ ]$\hspace{0.38cm} (A\ot \mu_{H})\co (\Psi\ot \mu_{H})\co (H\ot \Psi\ot H)\co (((H\ot \lambda_{H})\co \delta_{H})\ot A\ot H)$

\item[ ]$=((\varphi_{A}\co (\mu_{H}\ot A))\ot (\mu_{H}\co (\phi_{H}\ot \mu_{H})\co (H\ot \Psi\ot H)))\co (\delta_{H\ot H\ot A}\ot H)\co   (((H\ot \lambda_{H})$
\item[ ]$\hspace{0.38cm}\co \delta_{H})\ot A\ot H)$ {\scriptsize ({\blue by (\ref{hp})})} 

\item[ ]$=((\varphi_{A}\co (\mu_{H}\ot H)\co (H\ot \lambda_{H}\ot H))\ot (\mu_H\co (\phi_H\ot \mu_H)\co (H\ot \Psi\ot H)\co  (((H\ot \lambda_H)\co\delta_H)\ot A\ot H)))$
\item[ ]$\hspace{0.38cm}\co (H\ot H\ot  c_{H,A}\ot A\ot H)\co (H\ot c_{H,H}\ot A\ot A\ot H)\co (((\delta_{H}\ot H)\co \delta_{H})\ot \delta_{A}\ot H)$ {\scriptsize  ({\blue by }}
\item[ ]$\hspace{0.38cm} $ {\scriptsize  {\blue coassociativity of $\delta_{H}$, naturality of $c$ and (\ref{anticm}) for $\lambda_{H}$})}

\item[ ]$= ((\varphi_{A}\co (\mu_{H}\ot H)\co (H\ot \lambda_{H}\ot H))\ot  \varepsilon_{H}\ot \varepsilon_{A}\ot H)\co (H\ot H\ot  c_{H,A}\ot A\ot H)$
\item[ ]$\hspace{0.38cm} \co (H\ot c_{H,H}\ot A\ot A\ot H)\co (((\delta_{H}\ot H)\co \delta_{H})\ot \delta_{A}\ot H)$ {\scriptsize  ({\blue by (\ref{d6})})} 

\item[ ]$=(\varphi_{A}\co ((id_{H}\ast \lambda_{H})\ot A))\ot H$  {\scriptsize  ({\blue  by naturality of $c$ and counit properties})}

\item[ ]$=\varepsilon_{H}\ot A\ot H$  {\scriptsize  ({\blue by (\ref{2-conv}) for $H$ and the condition of left $H$-module for $A$}),}

\end{itemize}
and thus (\ref{adl2}) holds. The proof for (\ref{adl3}) is the following:

\begin{itemize}
\item[ ]$\hspace{0.38cm} (\mu_{A}\ot H)\co (\mu_{A}\ot \Psi)\co (A\ot \Psi\ot A)\co (A\ot H\ot (( \lambda_{A}\ot A)\co \delta_{A}))$

\item[ ]$=
((\mu_{A}\co (\mu_{A}\ot \varphi_{A})\co (H\ot \Psi\ot A))\ot (\phi_{H}\co (H\ot \mu_{A})))\co (A\ot \delta_{H\ot A\ot A})\co (A\ot H\ot ((\lambda_{A}\ot A)\co \delta_{A}))$
\item[ ]$\hspace{0.38cm}${\scriptsize ({\blue by (\ref{ap})})} 

\item[ ]$= ((\mu_{A}\co (\mu_{A}\ot \varphi_{A})\co (A\ot \Psi\ot A)\co (A\ot H\ot ((\lambda_{A}\ot A)\co \delta_{A})))\ot (\phi_{H}\co (H\ot \mu_{A})))$
\item[ ]$\hspace{0.38cm}\co (A\ot H\ot c_{H,A}\ot A\ot A)\co (A\ot H\ot H\ot c_{A,A}\ot A)\co (A\ot \delta_{H}\ot (((\lambda_{A}\ot A)\co \delta_{A})\co \delta_{A})$
\item[ ]$\hspace{0.38cm} $ {\scriptsize  ({\blue by coassociativity of $\delta_{A}$, naturality of $c$  and (\ref{anticm}) for $\lambda_{A}$})}

\item[ ]$=(A\ot \varepsilon_{H}\ot \varepsilon_{A}\ot (\phi_{H}\co (H\ot \mu_{A})))\co (A\ot H\ot c_{H,A}\ot A\ot A) \co (A\ot H\ot H\ot c_{A,A}\ot A)$
\item[ ]$\hspace{0.38cm} \co (A\ot \delta_{H}\ot (((\lambda_{A}\ot A)\co \delta_{A})\co \delta_{A})$ {\scriptsize  ({\blue by (\ref{d9})})} 

\item[ ]$=A\ot  (\phi_{H}\co (H\ot (\lambda_{A}\ast id_{A})))$  {\scriptsize  ({\blue  by naturality of $c$ and counit properties})}

\item[ ]$=A\ot H\ot \varepsilon_{A}$  {\scriptsize  ({\blue by (\ref{1-conv}) for $A$ and the condition of right $A$-module for $H$}),}

\end{itemize}

Finally, 
\begin{itemize}
\item[ ]$\hspace{0.38cm} (\mu_{A}\ot H)\co (\mu_{A}\ot \Psi)\co (A\ot \Psi\ot A)\co (A\ot H\ot (( A\ot \lambda_{A})\co \delta_{A}))$

\item[ ]$=((\mu_{A}\co (\mu_{A}\ot \varphi_{A})\co (H\ot \Psi\ot A))\ot (\phi_{H}\co (H\ot \mu_{A})))\co (A\ot \delta_{H\ot A\ot A})\co (A\ot H\ot ((A\ot \lambda_{A})\co \delta_{A}))$
\item[ ]$\hspace{0.38cm}$ {\scriptsize ({\blue by (\ref{ap})})} 

\item[ ]$=((\mu_{A}\co (\mu_{A}\ot \varphi_{A})\co (H\ot \Psi\ot A))\ot (\phi_{H}\co (H\ot (id_{A}\ast \lambda_{A})))\co (A\ot H\ot A\ot c_{H,A}\ot A) $  
\item[ ]$\hspace{0.38cm} \co (A\ot  H\ot c_{H,A}\ot c_{A,A})\co (A\ot \delta_{H}\ot ((\delta_{A}\ot \lambda_{A})\co \delta_{A}))$ {\scriptsize  ({\blue by coassociativity of $\delta_{A}$, naturality of $c$ )}}
\item[ ]$\hspace{0.38cm}$ {\scriptsize  {\blue and (\ref{anticm} for $\lambda_{A}$})}

\item[ ]$=((\mu_{A}\co (\mu_{A}\ot \varphi_{A})\co (H\ot \Psi\ot A))\ot H)\co (A\ot H\ot A\ot c_{H,A})\co (A\ot H\ot c_{H,A}\ot A)$
\item[ ]$\hspace{0.38cm} \co (A\ot \delta_{H}\ot ((A\ot \lambda_{A})\co \delta_{A}))$ {\scriptsize  ({\blue by (\ref{2-conv}) for $A$, counit properties, naturality of $c$ and the condition of left }}
\item[ ]$\hspace{0.38cm}$ {\scriptsize  {\blue $A$-module  for $H$})} 

\item[ ]$= ((\mu_A\co(\mu_{A}\ot \varphi_A)\co (A\ot \Psi\ot \lambda_{A})\co (A\ot H\ot \delta_{A}))\ot H)\co (A\ot H\ot c_{H,A})\co (A\ot \delta_{H}\ot A)$ 
\item[ ]$\hspace{0.38cm}$ {\scriptsize  ({\blue  by naturality of $c$})}

\item[ ]$=(A\ot \varepsilon_{H}\ot \varepsilon_{A}\ot H)\co (A\ot H\ot c_{H,A})\co (A\ot \delta_{H}\ot A) $ {\scriptsize  ({\blue by (\ref{d8})})}

\item[ ]$=A\ot H\ot \varepsilon_{A} $  {\scriptsize  ({\blue by naturality of $c$ and counit properties})}

\end{itemize}
and we obtain that (\ref{adl4}) holds. Therefore, $\Psi$ is  an $a$-comonoidal distributive law. 
}
\end{example}

\begin{example}
\label{2}
{\rm In this example, following Sections 4 and 5 of \cite{our2}, we provide examples of $a$-comonoidal distributive laws associated to skew parings between Hopf quasigroups. If $A$, $H$ are Hopf quasigroups, a skew pairing between $A$ and $H$ over $K$ is a morphism $\tau:A\ot H\rightarrow K$ such that the equalities
\begin{equation}
\label{skw1}\tau\co (\mu_A\ot H)=(\tau\ot \tau)\co (A\ot c_{A,H}\ot H)\co (A\ot A\ot\delta_H),
\end{equation}
\begin{equation}
\label{skw2}\tau\co (A\ot \mu_H)=(\tau\ot \tau)\co (A\ot c_{A,H}\ot H)\co ((c_{A,A}\co\delta_{A})\ot H\ot H),
\end{equation}
\begin{equation}
\label{skw3}\tau\co (A\ot \eta_{H})=\varepsilon_{A}, 
\end{equation}
\begin{equation}
\label{skw4}\tau\co (\eta_{A}\ot H)=\varepsilon_{H}, 
\end{equation}
hold.

If $\tau:A\ot H\rightarrow K$ is a skew pairing, by \cite[Proposition 4.3]{our2}, we have that $\tau$ is convolution invertible and, if $\tau^{-1}$ is the inverse of $\tau$, $\tau^{-1}=\tau\co (\lambda_{A}\ot H)$, $\tau=\tau^{-1}\co (A\ot \lambda_{H})$ and satisfies 
\begin{equation}
\label{skw5}
\tau^{-1}\co (A\ot \mu_{H})=(\tau^{-1}\ot \tau^{-1})\co (A\ot c_{A,H}\ot H)\co (\delta_{A}\ot H\ot H),
\end{equation}
\begin{equation}
\label{skw6}
\tau^{-1}\co (\mu_A\ot H)=(\tau^{-1}\ot \tau^{-1})\co (A\ot c_{A,H}\ot H)\co (A\ot A\ot (c_{H,H}\co \delta_H)), 
\end{equation}
(\ref{skw3}) and (\ref{skw4}). Therefore, $\tau: A\ot H\rightarrow K$ is a skew pairing between $A$ and $H$ if and only if $\tau^{-1}$ is a Hopf pairing between $A$ and $H$ in the sense of \cite{FT}. As a consequence of the properties of $\tau$ 
$$A\bowtie_{\tau} H=(A\ot H, \eta_{A\bowtie_{\tau} H}, \mu_{A\bowtie_{\tau} H}, \varepsilon_{A\bowtie_{\tau} H}, \delta_{A\bowtie_{\tau} H}, \lambda_{A\bowtie_{\tau} H})$$ is a Hopf quasigroup (see \cite[Proposition 2.2]{FT} and \cite[Corollary 4.10]{our2}) with unit, product, counit, coproduct and antipode defined by 
$$\eta_{A\bowtie_{\tau} H}=\eta_{A\otimes H},\;\; \mu_{A\bowtie_{\tau} H}=(\mu_A\ot\mu_H)\co (A\ot \Psi\ot H),$$
$$\varepsilon_{A\bowtie_{\tau} H}=\varepsilon_{A\otimes H},\;\; \delta_{A\bowtie_{\tau} H}=\delta_{A\otimes H},$$
$$\lambda_{A\bowtie_{\tau} H}=(\lambda_A\ot \lambda_H)\co \Psi\co c_{A,H}, 
$$
where 
\begin{equation}
\label{skwpsi}
\Psi=
(\tau\ot A\ot H\ot \tau^{-1})\co (A\ot H\ot \delta_{A\ot H})\co \delta_{A\ot H}\co c_{H,A}. 
\end{equation}

By \cite[Proposition 5.2]{our2}) we know that, if $\varphi_A:H\ot A\rightarrow A$ and $\phi_H:H\ot A\rightarrow H$ are defined as $$\varphi_A=(\tau\ot A\ot \tau^{-1})\co (A\ot H\ot \delta_A\ot H)\co \delta_{A\ot H}\co c_{H,A}$$ and $$\phi_H=(\tau\ot H\ot \tau^{-1})\co (A\ot H\ot c_{A,H}\ot H)\co (A\ot H\ot A\ot \delta_H)\co \delta_{A\ot H}\co c_{H,A}, $$
the pair $(A, \varphi_A)$ is a left $H$-module comonoid,  $(H, \phi_H)$ is a right $A$-module comonoid and 
$$\Psi=(\varphi_A\ot \phi_H)\co \delta_{H\ot A}.$$

As a consequence,  $A\bowtie_{\tau} H$ (or,  following the notation of \cite{FT}, $D(A,H,\tau^{-1})$) is the double cross product induced by the actions $\varphi_{A}$ and $\phi_{H}$. Therefore, by the previous example, the morphism $\Psi$ defined in (\ref{skwpsi}) is an $a$-comonoidal distributive law of $A$ over $H$. 
}
\end{example}

\begin{example}
\label{3}
{\rm This example is a particular case of the previous one. Let ${\Bbb F}$ be a field such that Char(${\Bbb F}$)$\neq 2$ and denote the tensor product over ${\Bbb F}$ as $\ot$. Consider the nonabelian group $S_{3}=\{\sigma_{0}, \sigma_{1}, \sigma_{2}, \sigma_{3}, \sigma_{4}, \sigma_{5}\}$, where $\sigma_{0}$ is  the identity, $o(\sigma_{1})=o(\sigma_{2})=o(\sigma_{3})=2$ and $o(\sigma_{4})=o(\sigma_{5})=3$. Let $u$ be an additional element such that $u^2=1$. By \cite[Theorem 1]{Chein}  the set 
$$L=M(S_{3},2)=\{\sigma_{i}u^{\alpha}\;; \; \alpha=0,1\}$$
is a Moufang loop where the product is defined by 
$$\sigma_{i}u^{\alpha}.\;\sigma_{j}u^{\beta}=(\sigma_{i}^{\nu}\sigma_{j}^{\mu})^{\nu}u^{\alpha +\beta},\;\;\nu=(-1)^{\beta}, \; \mu=(-1)^{\alpha +\beta}.$$

Then, $L$ is an I.P. loop and  ${\Bbb F}L$ is a cocommutative Hopf quasigroup. 
 On the other hand, let $H_{4}$ be the $4$-dimensional Taft Hopf algebra.  The basis of $H_{4}$ is $\{1,x,y,w=xy\}$ and the multiplication table is defined by 
\begin{center}
\begin{tabular}{|c|c|c|c|c|}
\hline  $\;$ & $x$ & $y$ & $w$   \\
\hline  $ x$ &  $1$ & $w$ & $y$ \\
\hline  $ y$ &  $-w$ & $0$ & $0$ \\
\hline  $ w$ &  $-y$ & $0$ & $0$ \\
\hline
\end{tabular}
\end{center}

The costructure of $H_{4}$ is given by 
$$ \delta_{H_{4}}(x)=x\ot x,\; \delta_{H_{4}}(y)=y\ot x +1\ot y,\; \delta_{H_{4}}(w)=w\ot 1 +x\ot w,$$
$$\varepsilon_{H_{4}}(x)=1_{\Bbb F},\; \varepsilon_{H_{4}}(y)=\varepsilon_{H_{4}}(w)=0$$
and the antipode $\lambda_{H_{4}}$ is described by 
$$\l \lambda_{H_{4}}(x)=x,\; \lambda_{H_{4}}(y)=w,\; \lambda_{H_{4}}(w)=-y.$$

The morphism $\tau:{\Bbb F}L\ot H_{4}\rightarrow {\Bbb F}$ defined by 
$$\tau (\sigma_{i}u^{\alpha}\ot z)=\left\{ \begin{array}{ccc} 1 & {\rm if} &
z=1 \\
 (-1)^{\alpha} & {\rm if} &
z=x \\ 
 0 & {\rm if} &
z=y, w 
 \end{array}\right.$$
 is a skew pairing such that $\tau=\tau^{-1}$. Then, ${\Bbb F}L\bowtie_{\tau} H_{4}$ is a Hopf quasigroup where the $a$-comonoidal distributive law $\Psi$ of ${\Bbb F}L$ over $H_{4}$ is defined by:
 $$\Psi(1\ot \sigma_{i}u^{\alpha})=\sigma_{i}u^{\alpha}\ot 1, \;\; \Psi(x\ot \sigma_{i}u^{\alpha})=\sigma_{i}u^{\alpha}\ot x,$$ $$ \Psi(y\ot \sigma_{i}u^{\alpha})=(-1)^{\alpha}\sigma_{i}u^{\alpha}\ot y,  \;\;\Psi(w\ot \sigma_{i}u^{\alpha})=(-1)^{\alpha}\sigma_{i}u^{\alpha}\ot w. $$
}
\end{example}

\begin{example}
\label{4}
{\rm  Let $A$, $H$ be Hopf quasigroups and assume that there exists an action $\varphi_{A}:H\ot A\rightarrow A$ such that $(A, \varphi_{A})$  is a left $H$-quasimodule magma, $(A, \varphi_{A})$ a left $H$-quasimodule comonoid and satisfies the identities
\begin{equation}
\label{cesp1}
(H\ot \varphi_{A})\co ((c_{H,H}\co \delta_{H})\ot A)=(H\ot \varphi_{A})\co (\delta_{H}\ot A),  
\end{equation}
\begin{equation}
\label{cesp2}
\varphi_{A}\co (H\ot (\varphi_{A}\co (\lambda_{H}\ot A)))=\varphi_{A}\co ((\mu_{H}\co (H\ot \lambda_{H}))\ot A).  
\end{equation}

Define $\Psi:H\otimes A\rightarrow A\otimes H$ by  
\begin{equation}
\label{R}
\Psi=(\varphi_{A}\ot H)\co (H\ot c_{H,A})\co (\delta_{H}\ot A).
\end{equation}

Then, by the monoidal version of the proofs  developed in \cite{BrzezJiao2}, we can assure that the morphism $\Psi$ satisfies (\ref{dl1-1}), 
\begin{equation}
\label{dl21}
\Psi\co (\mu_{H}\ot A)\co (H\ot \lambda_{H}\ot A)=(A\ot \mu_{H})\co ( \Psi\ot H)\co (H\ot \Psi)\co (H\ot \lambda_{H}\ot A),
\end{equation}
(\ref{dl3}) and (\ref{dl4}). Therefore, is an example of distributive law of $H$ over $A$. Moreover, $\Psi$ is comonoidal and we have that 
\begin{equation}
\label{l-con}
(\varepsilon_{A}\ot H)\co \Psi=H\ot \varepsilon_{A}.
\end{equation}

Also,  the proof of \cite[Lemma 2.6]{BrzezJiao2} remains valid  in this general symmetric monoidal setting and, as a consequence of the properties of $\Psi$, the equality
\begin{equation}
\label{cesp3}
\Psi\circ (A\ot \lambda_{A})=(\lambda_{A}\ot H)\co \Psi
\end{equation}
holds. 

Finally, $\Psi$ is a $a$-comonoidal distributive law because:

\begin{itemize}
\item[ ]
\item[ ]$\hspace{0.38cm} (A\ot \mu_{H})\co (\Psi\ot \mu_{H})\co (H\ot \Psi\ot H)\co (((\lambda_{H}\ot H)\co \delta_{H})\ot A\ot H)$

\item[ ]$=((\varphi_{A}\co (\lambda_{H}\ot \varphi_{A})\co (\delta_{H}\ot A))\ot (\mu_H\circ (\lambda_{H}\ot \mu_H)))\co (H\ot  c_{H,A}\ot H\ot H)\co  (c_{H,H}\ot c_{H,A}\ot  H) $ 
\item[ ]$\hspace{0.38cm}\co (((\delta_{H}\ot H)\co \delta_{H})\ot A\ot H)$ {\scriptsize ({\blue by naturality of $c$, coassociativity of $\delta_{H}$ and (\ref{anticm}) for $\lambda_{H}$})}

\item[ ]$=(A\ot (\mu_H\circ (\lambda_{H}\ot \mu_H)\circ (\delta_H\ot H)))\co (c_{H,A}\ot H)$  {\scriptsize  ({\blue by (\ref{pq}) for $\varphi_{A}$})}

\item[ ]$=\varepsilon_{H}\ot A\ot H$  {\scriptsize  ({\blue by (\ref{lH})  and naturality of $c$}),}

\item[ ]

\end{itemize}

\begin{itemize}
\item[ ]$\hspace{0.38cm} (A\ot \mu_{H})\co (\Psi\ot \mu_{H})\co (H\ot \Psi\ot H)\co (((H\ot \lambda_{H})\co \delta_{H})\ot A\ot H)$

\item[ ]$=((\varphi_{A}\co (H\ot (\varphi_{A}\co (\lambda_{H}\ot A))))\ot (\mu_H\circ (H\ot \mu_H)\circ (H\ot \lambda_H\ot H)\circ (\delta_H\ot H)))$ 
\item[ ]$\hspace{0.38cm}\co (H\ot H\ot c_{H,A}\ot H)\co (H\ot c_{H,H}\ot A\ot H)\co (((\delta_{H}\ot H)\co \delta_{H})\ot A\ot H)$ {\scriptsize ({\blue by naturality of }} 
\item[ ]$\hspace{0.38cm}${\scriptsize {\blue $c$, coassociativity of $\delta_{H}$ and (\ref{anticm}) for $\lambda_{H}$})} 

\item[ ]$=(\varphi_{A}\co (H\ot \varphi_{A})\co (((H\ot \lambda_{H})\co \delta_{H})\ot A)\ot H$  {\scriptsize  ({\blue by (\ref{lH}), naturality of $c$ and counit properties})}

\item[ ]$=\varepsilon_{H}\ot A\ot H$ {\scriptsize  ({\blue by (\ref{pq}) for $\varphi_{A}$}),}

\item[ ]

\end{itemize}

\begin{itemize}
\item[ ]$\hspace{0.38cm}  (\mu_{A}\ot H)\co (\mu_{A}\ot \Psi)\co (A\ot \Psi\ot A)\co (A\ot H\ot ((  \lambda_{A}\ot A)\co \delta_{A}))$

\item[ ]$= (((\mu_{A}\co ((\mu_{A}\co (A\ot \lambda_{A}))\ot A))\co (A\ot ((\varphi_{A}\ot \varphi_{A})\co \delta_{H\ot A})))\ot H)\co (A\ot H\ot c_{H,A})\co (A\ot \delta_{H}\ot A)$ 
\item[ ]$\hspace{0.38cm}$ {\scriptsize ({\blue by (\ref{cesp3}), naturality of $c$ and coassociativity of $\delta_{H}$})}

\item[ ]$=(((\mu_{A}\co ((\mu_{A}\co (A\ot \lambda_{A}))\ot A))\co (A\ot (\delta_{A}\co \varphi_{A})))\ot H)\co (A\ot H\ot c_{H,A})\co (A\ot \delta_{H}\ot A)$
\item[ ]$\hspace{0.38cm}$ {\scriptsize ({\blue by (\ref{dq})})} 

\item[ ]$=A\ot (((\varepsilon_{A}\co \varphi_{A})\ot H)\co (H\ot c_{H,A})\co (\delta_{H}\ot A))$ {\scriptsize ({\blue by (\ref{rH})})} 

\item[ ]$=A\ot H\ot \varepsilon_{A}$  {\scriptsize  ({\blue by (\ref{eq}), naturality of $c$ and counit properties}),}

\item[ ]

\end{itemize}

\begin{itemize}
\item[ ]$\hspace{0.38cm}  (\mu_{A}\ot H)\co (\mu_{A}\ot \Psi)\co (A\ot \Psi\ot A)\co (A\ot H\ot (( A\ot \lambda_{A})\co \delta_{A}))$

\item[ ]$=(((\mu_{A}\co (\mu_{A}\ot \lambda_{A}))\co (A\ot ((\varphi_{A}\ot \varphi_{A})\co \delta_{H\ot A})))\ot H)\co (A\ot H\ot c_{H,A})\co (A\ot \delta_{H}\ot A)$ 
\item[ ]$\hspace{0.38cm}$ {\scriptsize ({\blue by (\ref{cesp3}), naturality of $c$ and coassociativity of $\delta_{H}$})} 

\item[ ]$=(((\mu_{A}\co (\mu_{A}\ot \lambda_{A}))\co (A\ot (\delta_{A}\co \varphi_{A})))\ot H)\co (A\ot H\ot c_{H,A})\co (A\ot \delta_{H}\ot A)$ {\scriptsize ({\blue by (\ref{dq})})}

\item[ ]$=A\ot (((\varepsilon_{A}\co \varphi_{A})\ot H)\co (H\ot c_{H,A})\co (\delta_{H}\ot A))$ {\scriptsize ({\blue by (\ref{rH})})} 

\item[ ]$=A\ot H\ot \varepsilon_{A}$  {\scriptsize  ({\blue by (\ref{eq}), naturality of $c$ and counit properties}).}

\end{itemize}

Let ${\Bbb F}$ be a field  and denote the tensor product over ${\Bbb F}$ as $\ot$. Let $A$, $H$ be Hopf quasigroups in ${\Bbb F}$-{\sf Vect}. The morphism $\Psi$, defined in (\ref{R}), was used in \cite{BrzezJiao1} to introduce the notion of  smash product for the quasigroups $H$ and $A$. This product is a particular instance of  a more general construction that we can find in \cite{BrzezJiao2}. In \cite{BrzezJiao2}, for a morphism $R:H\ot A\rightarrow A\ot H$,  the authors introduce  a product in $A\ot H$, defined by 
$$\mu=(\mu_{A}\ot \mu_{H})\co (A\ot R\ot H),$$
and an antipode $\lambda=R\co (\lambda_{H}\ot \lambda_{A})\co c_{A,H}$. With this new product and the antipode $\lambda$, they proved that, with the tensor product coproduct and the tensor product counit, if $R$ satisfies (\ref{dl1-1}) and (\ref{l-con}), $A\ot H$ has a new structure of Hopf quasigroup (called the $R$-smash product Hopf quasigroup for $H$ and $A$) if and only if $R$ is a comonoid morphism, satisfies (\ref{dl3}), (\ref{dl4}), (\ref{dl21}) and 
$$(A\ot \varepsilon_{H})\co R\co c_{A,H}\co (A\ot \lambda_{H})\co R=\varepsilon_{H}\ot A.$$

Note that the condition (\ref{l-con}) in the Example \ref{1} implies that 
$$\phi_{H}=H\ot \varepsilon_{A},$$
i.e., $\phi_{H}$ is trivial. Then, the theory of $R$-smash products presented in \cite{BrzezJiao2} does not cover as a particular cases the  double cross products $A\bowtie H$ defined in \cite{our2}.

}
\end{example}

\begin{example}
\label{5}
{\rm Let $A$, $H$ be Hopf quasigroups and assume that there exists two actions $\varphi_{A}:H\ot A\rightarrow A$, $\hat{\varphi}_{A}:H\ot A\rightarrow A$ such that $(A, \varphi_{A})$,   $(A, \hat{\varphi}_{A})$ are  left $H$-quasimodule magmas and  left $H$-quasimodule comonoids and satisfy the identities (\ref{cesp1}),  (\ref{cesp2}) and
\begin{equation}
\label{cesp4}
\hat{\varphi}_{A}\co (H\ot \varphi_{A})=\varphi_{A}\co (H\ot \hat{\varphi}_{A})\co (c_{H,H}\ot A).
\end{equation}

Define $\Gamma:H\otimes A\rightarrow A\otimes H$ by  
\begin{equation}
\label{R1}
\Gamma=(\varphi_{A}\ot H)\co (H\ot T)\co (\delta_{H}\ot A), 
\end{equation}
where 
\begin{equation}
\label{R2}
T=(\hat{\varphi}_{A}\ot H)\co (H\ot c_{H,A})\co (\delta_{H}\ot A).
\end{equation}

By the previous example, we known that the morphisms $\Psi$ (defined in (\ref{R})) and $T$ are examples of $a$-comonoidal distributive laws. Moreover, they satisfy the identities (\ref{dl1-1}), (\ref{dl21}), (\ref{l-con}) and (\ref{cesp3}). On the other hand, by the naturality of $c$ and the coassociativity of $\delta_{H}$, we obtain 
\begin{equation}
\label{R3}
(A\ot \delta_{H})\co T=(T\ot H)\co (H\ot c_{H,A})\co (\delta_{H}\ot A)
\end{equation}
and, using (\ref{cesp4}), again the naturality of $c$ and the coassociativity of $\delta_{H}$,   we can prove that 
\begin{equation}
\label{R4}
 T\co (H\ot \varphi_{A})=(\varphi_{A}\ot H)\co (H\ot T)\co (c_{H,H}\ot A)
\end{equation}
holds. Then, 
\begin{equation}
\label{R5}
 (c_{H,A}\ot H)\co (H\ot T)\co (\delta_{H}\ot A)= (T\ot H)\co (H\ot c_{H,A})\co (\delta_{H}\ot A).
\end{equation}
Indeed:
\begin{itemize}
\item[ ]$\hspace{0.38cm}  (c_{H,A}\ot H)\co (H\ot T)\co (\delta_{H}\ot A)$

\item[ ]$=(c_{H,A}\ot H)\co (H\ot \hat{\varphi}_{A}\ot H)\co (\delta_{H}\ot c_{H,A})\co (\delta_{H}\ot A)$ {\scriptsize ({\blue by coassociativity of $\delta_{H}$})} 

\item[ ]$=(c_{H,A}\ot H)\co (H\ot \hat{\varphi}_{A}\ot H)\co ((c_{H,H}\co\delta_{H})\ot c_{H,A})\co (\delta_{H}\ot A) $ {\scriptsize   ({\blue by  (\ref{cesp1})})}

\item[ ]$= (T\ot H)\co (H\ot c_{H,A})\co (\delta_{H}\ot A)$  {\scriptsize  ({\blue by naturality of $c$}).}

\end{itemize}

Also, $T$ satisfies 
\begin{equation}
\label{R6}
(\delta_{A}\ot H)\co T= (A\ot T)\co (T\ot A)\co (H\ot \delta_{A})
\end{equation}
because
\begin{itemize}
\item[ ]$\hspace{0.38cm}  (\delta_{A}\ot H)\co T$

\item[ ]$=(((\hat{\varphi}_{A}\ot \hat{\varphi}_{A})\co \delta_{H\ot A})\ot H)\co (H\ot c_{H,A})\co (\delta_{H}\ot A) $
{\scriptsize ({\blue by  (\ref{dq}) for $\hat{\varphi}_{A}$)})} 

\item[ ]$= (A\ot T)\co (T\ot A)\co (H\ot \delta_{A})$  {\scriptsize  ({\blue by naturality of $c$ and coassociativity of $\delta_{H}$})}

\end{itemize}
and, as a consequence of this fact, 
\begin{equation}
\label{R7}
(\delta_{A}\ot H)\co \Gamma= (A\ot \Gamma)\co (\Gamma\ot A)\co (H\ot \delta_{A})
\end{equation}
also holds. Indeed:
\begin{itemize}
\item[ ]$\hspace{0.38cm}  (\delta_{A}\ot H)\co \Gamma$

\item[ ]$=(((\varphi_{A}\ot \varphi_{A})\co \delta_{H\ot A})\ot H)\co (H\ot T)\co (\delta_{H}\ot A) $
{\scriptsize ({\blue by  (\ref{dq}) for $\varphi_{A}$)})} 

\item[ ]$= (\varphi_{A}\ot\varphi_{A}\ot H)\co (H\ot c_{H,A}\ot T)\co (\delta_{H}\ot T\ot A) \co (\delta_{H}\ot \delta_{A})$  {\scriptsize  ({\blue by (\ref{R6})})}

\item[ ]$= (\varphi_{A}\ot\varphi_{A}\ot H)\co (H\ot ((c_{H,A}\ot H)\co (H\ot T)\co (\delta_{H}\ot A))\ot A) \co (\delta_{H}\ot \delta_{A})$  {\scriptsize  ({\blue by  coassociativity}}
\item[ ]$\hspace{0.38cm}$ {\scriptsize  {\blue of $\delta_{H}$})}

\item[ ]$= (A\ot \Gamma)\co (\Gamma\ot A)\co (H\ot \delta_{A})$  {\scriptsize  ({\blue by (\ref{R3}) and (\ref{R5})}).}

\end{itemize}

The morphism $\Gamma $ is a distributive law of $H$ over $H$ because satisfy the identities (\ref{dl1-1}), (\ref{dl21}), (\ref{dl3}) and (\ref{dl4}). Indeed,  (\ref{dl1-1}) holds because 

\begin{itemize}
\item[ ]$\hspace{0.38cm}(\mu_{A}\ot H)\co (A\ot \Gamma)\co (\Gamma\ot A)  $

\item[ ]$= ((\mu_{A}\co (\varphi_{A}\ot \varphi_{A}))\ot H)\co (H\ot A\ot H\ot T)\co (H\ot ((A\ot \delta_{H})\co T)\ot A)\co (\delta_{H}\ot A\ot A)$ 
\item[ ]$\hspace{0.38cm}$ {\scriptsize ({\blue by definition})} 

\item[ ]$= ((\mu_{A}\co (\varphi_{A}\ot \varphi_{A}))\ot H)\co (H\ot A\ot H\ot T)\co (H\ot ((T\ot H)\co (H\ot c_{H,A})\co (\delta_{H}\ot A))\ot A)$
\item[ ]$\hspace{0.38cm}\co (\delta_{H}\ot A\ot A)$ {\scriptsize ({\blue by  (\ref{R3})})} 

\item[ ]$= ((\mu_{A}\co (\varphi_{A}\ot \varphi_{A}))\ot H)\co (H\ot A\ot H\ot T)\co (H\ot ( (c_{H,A}\ot H)\co (H\ot T)\co (\delta_{H}\ot A))\ot A)$
\item[ ]$\hspace{0.38cm}\co (\delta_{H}\ot A\ot A)$ {\scriptsize ({\blue by  (\ref{R5})})}

\item[ ]$=((\mu_A\co \varphi_{A\ot A})\ot H)\co (H\ot A\ot T)\co (H\ot T\ot A)\co (\delta_{H}\ot A\ot A)$  {\scriptsize  ({\blue by coassociativity of $\delta_{H}$})}

\item[ ]$=((\varphi_A\co (H\ot \mu_A))\ot H)\co (H\ot A\ot T)\co (H\ot T\ot A)\co (\delta_{H}\ot A\ot A) $ {\scriptsize ({\blue by (\ref{muq})})} 

\item[ ]$= \Gamma\co (H\ot \mu_{A})${\scriptsize ({\blue by  (\ref{dl1-1}) for $T$.})} 

\end{itemize}

On the other hand, 
\begin{itemize}
\item[ ]$\hspace{0.38cm} (A\ot \mu_{H})\co ( \Gamma\ot H)\co (H\ot \Gamma)\co (H\ot \lambda_{H}\ot A) $

\item[ ]$=(\varphi_{A}\ot \mu_{H})\co (H\ot (T\co (H\ot \varphi_{A}))\ot H)\co (H\ot H\ot H\ot T)\co (\delta_{H}\ot (\delta_{H}\co \lambda_{H})\ot A) $ 
\item[ ]$\hspace{0.38cm}${\scriptsize ({\blue by definition})} 

\item[ ]$= (\varphi_{A}\ot \mu_{H})\co (H\ot ((\varphi_{A}\ot H)\co (H\ot T)\co (c_{H,H}\ot A))\ot H)\co (H\ot H\ot H\ot T)$
\item[ ]$\hspace{0.38cm}\co (\delta_{H}\ot (\delta_{H}\co \lambda_{H})\ot A)${\scriptsize ({\blue by  (\ref{R4})})} 

\item[ ]$= ((\varphi_{A}\co (H\ot (\varphi_{A}\co (\lambda_{H}\ot A))))\ot H)\co (H\ot H\ot ( (A\ot \mu_{H})\co ( T\ot H)\co (H\ot T)\co (H\ot \lambda_{H}\ot A) ))$
\item[ ]$\hspace{0.38cm}\co (H\ot c_{H,H}\ot H\ot A)\co (\delta_{H}\ot (c_{H,H}\co \delta_{H})\ot A)$ {\scriptsize ({\blue by  naturality of $c$ and (\ref{anticm}) for $\lambda_{H}$})}

\item[ ]$=((\varphi_{A}\co ((\mu_{H}\co (H\ot \lambda_{H}))\ot A))\ot H)\co (H\ot H\ot ( T\co ((\mu_{H}\co (H\ot \lambda_{H}))\ot A))) $  
\item[ ]$\hspace{0.38cm}\co  (H\ot c_{H,H}\ot H\ot A)\co (\delta_{H}\ot (c_{H,H}\co \delta_{H})\ot A)$ {\scriptsize  ({\blue (\ref{cesp2}) for $\varphi_{A}$ and (\ref{dl21}) for $T$})}

\item[ ]$=(\varphi_{A}\ot H)\co (H\ot T)\co (((\mu_{H}\ot \mu_{H})\co \delta_{H\ot H}\co (H\ot \lambda_{H}))\ot A) $ {\scriptsize ({\blue by naturality of $c$ and (\ref{anticm}) for }} 
 \item[ ]$\hspace{0.38cm}$ {\scriptsize {\blue $\lambda_{H}$})} 

\item[ ]$= \Gamma\co (\mu_{H}\ot A)\co (H\ot \lambda_{H}\ot A) ${\scriptsize ({\blue by  (\ref{delta-mu})})} 

\end{itemize}
and, as a consequence, (\ref{dl21}) holds. The identities (\ref{dl3}) and (\ref{dl4}) for $\Gamma$ follow from (\ref{dl3}) and (\ref{dl4}) for $T$,  (\ref{etaq}),  (\ref{uq}),  (\ref{delta-eta}), the naturality of $c$ and the counit properties.

It is easy to show that $\varepsilon_{A\ot H}\co \Gamma=\varepsilon_{H\ot A}$ follows by (\ref{eq}), (\ref{cdl2}) for $T$ and the properties of the counit. Also, 

\begin{itemize}
\item[ ]$\hspace{0.38cm}\delta_{A\ot H}\ot \Gamma  $

\item[ ]$=(A\ot  c_{A,H}\ot H)\co (((\varphi_{A}\ot \varphi_{A})\co \delta_{H\ot A})\ot \delta_{H})\co (H\ot T)\co (\delta_{H}\ot A)$ 
{\scriptsize ({\blue by  (\ref{dq})})} 

\item[ ]$= (A\ot  c_{A,H}\ot H)\co (((\varphi_{A}\ot \varphi_{A})\co (H\ot c_{H,A}\ot A))\ot \delta_{H})\co (\delta_{H}\ot A\ot T)\co (H\ot T\ot A)\co (\delta_{H}\ot \delta_{A})$ 
\item[ ]$\hspace{0.38cm}${\scriptsize ({\blue by  (\ref{R6})})} 

\item[ ]$= (A\ot H\ot \varphi_{A}\ot H)\co (\varphi_{A}\ot c_{H,H}\ot A\ot H)\co (H\ot c_{H,A}\ot H\ot A\ot H)\co (\delta_{H}\ot A\ot H\ot T)$
\item[ ]$\hspace{0.38cm}\co (H\ot ((c_{H,A}\ot H)\co (H\ot T)\co (\delta_{H}\ot A))\ot A)\co (\delta_{H}\ot \delta_{A}) ${\scriptsize ({\blue by  (\ref{R3}) and  (\ref{R5})})} 

\item[ ]$=(\varphi_{A}\ot ((H\ot \varphi_{A})\co ((c_{H,H}\co \delta_{H})\ot A))\ot H)\co (H\ot ((c_{H,A}\ot T)\co (H\ot T\ot A)\co (\delta_{H}\ot \delta_{A}))) $
\item[ ]$\hspace{0.38cm}\co (\delta_{H}\ot A)  ${\scriptsize ({\blue by  naturality of $c$})}

\item[ ]$=(\varphi_{A}\ot ((H\ot \varphi_{A})\co (\delta_{H}\ot A))\ot H)\co (H\ot ((T\ot T)\co \delta_{H\ot A})))\co (\delta_{H}\ot A) $  
 {\scriptsize  ({\blue (\ref{cesp1})})}

\item[ ]$=(\varphi_{A}\ot H\ot \varphi_{A}\ot H)\co (H\ot T\ot H\ot T)\co (\delta_{H}\ot ((c_{H,A}\ot H)\co (H\ot c_{H,A})\co (\delta_{H}\ot A))\ot A) $ 
\item[ ]$\hspace{0.38cm}\co (\delta_{H}\ot \delta_{A})$ {\scriptsize ({\blue by  (\ref{R3}) and coassociativity of $\delta_{H}$})} 

\item[ ]$= (\Gamma\ot \Gamma)\co \delta_{H\ot A} $ {\scriptsize ({\blue by  naturality of $c$}),} 

\end{itemize}
and this implies that $\Gamma $ is comonoidal. 

The morphism $\Gamma$ also satisfies (\ref{l-con}) because (\ref{eq}) holds and $T$ satisfies (\ref{l-con}). Moreover, as in the previous example, we can assure that (\ref{cesp3}) holds for $\Gamma$. Then, 

\begin{itemize}
\item[ ]
\item[ ]$\hspace{0.38cm} (A\ot \mu_{H})\co (\Gamma\ot \mu_{H})\co (H\ot \Gamma\ot H)\co (((\lambda_{H}\ot H)\co \delta_{H})\ot A\ot H)$

\item[ ]$= (\varphi_{A}\ot \mu_{H})\co (H\ot (T\co (H\ot \varphi_{A}))\ot \mu_{H})\co ((\delta_{H}\co \lambda_{H})\ot H\ot T\ot H)\co (((H\ot \delta_{H})\co \delta_{H})\ot A\ot H)$ 
\item[ ]$\hspace{0.38cm}${\scriptsize ({\blue by definition})}

\item[ ]$=(\varphi_{A}\ot \mu_{H})\co (H\ot ((\varphi_{A}\ot H)\co (H\ot T)\co (c_{H,H}\ot A))\ot \mu_{H})\co ((\delta_{H}\co \lambda_{H})\ot H\ot T\ot H)$
\item[ ]$\hspace{0.38cm} \co (((H\ot \delta_{H})\co \delta_{H})\ot A\ot H)$  {\scriptsize  ({\blue by (\ref{R4})})}

\item[ ]$= ((\varphi_{A}\co (H\ot\varphi_{A})\co  (((\lambda_{H}\ot H)\co \delta_H)\ot A))\ot \mu_{H})\co (H\ot T\ot \mu_{H})\co (((H\ot \lambda_{H})\co c_{H,H}$
\item[ ]$\hspace{0.38cm}\co \delta_{H})\ot T\ot H) \co (\delta_{H}\ot A\ot H) $  {\scriptsize  ({\blue by (\ref{anticm}), coassociativity of $\delta_{H}$ and  naturality of $c$})}

\item[ ]$=(A\ot \mu_{H})\co (T\ot \mu_{H})\co (H\ot T\ot H)\co (((\lambda_{H}\ot H)\co \delta_{H})\ot A\ot H) $  {\scriptsize  ({\blue by (\ref{pq}) for $\varphi_{A}$, naturality}}
\item[ ]$\hspace{0.38cm}${\scriptsize  {\blue  of $c$ and counit properties})}

\item[ ]$=\varepsilon_{H}\ot A\ot H$  {\scriptsize  ({\blue by (\ref{adl1}) for T}),}

\item[ ]

\end{itemize}

\begin{itemize}
\item[ ]$\hspace{0.38cm} (A\ot \mu_{H})\co (\Gamma\ot \mu_{H})\co (H\ot \Gamma\ot H)\co (((H\ot \lambda_{H})\co \delta_{H})\ot A\ot H)$

\item[ ]$=(\varphi_{A}\ot \mu_{H})\co (H\ot (T\co (H\ot \varphi_{A}))\ot \mu_{H})\co (\delta_{H}\ot H\ot T\ot H)\co (((H\ot (\delta_{H}\co \lambda_{H}))\co \delta_{H})\ot A\ot H)$
\item[ ]$\hspace{0.38cm}${\scriptsize ({\blue by definition})}

\item[ ]$=(\varphi_{A}\ot \mu_{H})\co (H\ot ((\varphi_{A}\ot H)\co (H\ot T)\co (c_{H,H}\ot A))\ot \mu_{H})\co (\delta_{H}\ot H\ot T\ot H)$
\item[ ]$\hspace{0.38cm}\co (((H\ot (\delta_{H}\co \lambda_{H}))\co \delta_{H})\ot A\ot H) $  {\scriptsize  ({\blue by (\ref{R4})})}

\item[ ]$=((\varphi_{A}\co (H\ot \varphi_{A}))\ot H)\co (H\ot H\ot  ((A\ot \mu_{H})\co (T\ot \mu_{H})\co (H\ot T\ot H)\co (((H\ot \lambda_{H})$
\item[ ]$\hspace{0.38cm}\co \delta_{H})\ot A\ot H))) \co (((H\ot c_{H,H})\co (\delta_{H}\ot \lambda_{H})\co \delta_{H})\ot A\ot H)$ {\scriptsize  ({\blue by (\ref{anticm}), coassociativity of $\delta_{H}$}}
\item[ ]$\hspace{0.38cm}${\scriptsize  {\blue  and  naturality of $c$})}

\item[ ]$=(\varphi_{A}\co (H\ot\varphi_{A})\co  ((( H\ot \lambda_{H})\co \delta_H)\ot A))\ot H$ {\scriptsize  ({\blue by (\ref{adl2}) for $T$, naturality of $c$ and counit properties})}

\item[ ]$=\varepsilon_{H}\ot A\ot H$  {\scriptsize  ({\blue by (\ref{pq}) for $\varphi_{A}$}),}

\item[ ]

\end{itemize}

\begin{itemize}
\item[ ]$\hspace{0.38cm}  (\mu_{A}\ot H)\co (\mu_{A}\ot \Gamma)\co (A\ot \Gamma\ot A)\co (A\ot H\ot ((  \lambda_{A}\ot A)\co \delta_{A}))$

\item[ ]$=  (\mu_{A}\ot H)\co ((\mu_{A}\co (A\ot \lambda_{A}))\ot \Gamma)\co (A\ot \Gamma\ot A)\co (A\ot H\ot \delta_{A})$ 
{\scriptsize ({\blue by (\ref{cesp3}) for $\Gamma$})}

\item[ ]$=((\mu_A\circ (\mu_A\ot A)\circ (A\ot \lambda_A\ot A)\circ (A\ot \delta_A))\ot H)\co (A\ot \Gamma) ${\scriptsize ({\blue by (\ref{R7})})} 

\item[ ]$= A\ot ((\varepsilon_{A}\ot H)\co \Gamma) $ {\scriptsize ({\blue by (\ref{rH}) for $\varphi_{A}$})} 

\item[ ]$=A\ot H\ot \varepsilon_{A}$  {\scriptsize  ({\blue by (\ref{l-con}) for $\Gamma$}),}

\item[ ]

\end{itemize}

\begin{itemize}
\item[ ]$\hspace{0.38cm}  (\mu_{A}\ot H)\co (\mu_{A}\ot \Gamma)\co (A\ot \Gamma\ot A)\co (A\ot H\ot (( A\ot \lambda_{A})\co \delta_{A}))$

\item[ ]$=((\mu_{A}\co (A\ot \lambda_{A}))\ot H)\co (\mu_{A}\ot \Gamma)\co (A\ot \Gamma\ot A)\co (A\ot H\ot \delta_{A}) $ 
{\scriptsize ({\blue by (\ref{cesp3}) for $\Gamma$})} 

\item[ ]$=((\mu_A\circ(\mu_A\ot \lambda_A)\circ (A\ot \delta_A))\ot H)\co (A\ot \Gamma) ${\scriptsize ({\blue by (\ref{R7})})} 

\item[ ]$=  A\ot ((\varepsilon_{A}\ot H)\co \Gamma) $ {\scriptsize ({\blue by (\ref{rH}) for $\varphi_{A}$})} 

\item[ ]$=A\ot H\ot \varepsilon_{A}$  {\scriptsize  ({\blue by (\ref{l-con}) for $\Gamma$}).}

\end{itemize}

Therefore, $\Gamma$ is a an $a$-comonoidal distributive law of $H$ over $A$.

}
\end{example}

\begin{example}
\label{6}
{\rm Let $(A, \varphi_{A})$ be a left $H$-quasimodule magma-comonoid and assume that there exists a right action $\phi_{A}:A\ot H\rightarrow A$ such that $(A, \phi_{A})$ is a right $H$-quasimodule magma-comonoid. Following \cite{FT}, we will say that $(A, \varphi_{A}, \phi_{A})$ is a $H$-biquasimodule Hopf quasigroup if 
\begin{equation}
\label{61}
\varphi_{A}\co (H\ot \phi_{A})=\phi_{A}\co (\varphi_{A}\ot H)
\end{equation}
holds.

The main result of \cite{FT1} asserts that, if (\ref{cesp2}) and 
\begin{equation}
\label{62}
\phi_{A}\co ((\phi_{A}\co (A\ot \lambda_{H}))\ot H)=\phi_{A}\co (A\ot (\mu_{H}\co (\lambda_{H}\ot H))
\end{equation}
hold, the twisted smash product $A\circledast H$ built on $A\ot H$ with tensor coproduct and unit, where the product and antipode  are defined by 
\begin{equation}
\label{63}
\mu_{A\circledast H}=
\end{equation}
$$ (\mu_{A}\ot \mu_{H})\co (A\ot ((\phi_{A}\ot H)\co (A\ot ((\lambda_{H}\ot H)\co c_{H,H}\co \delta_{H}))\co 
(\varphi_{A}\ot H)\co (H\ot c_{H,A})\co (\delta_{H}\ot A))\ot H)$$
and 
\begin{equation}
\label{64}
\lambda_{A\circledast H}=
\end{equation}
$$(\phi_{A}\ot H)\co (A\ot ((\lambda_{H}\ot H)\co c_{H,H}\co \delta_{H}))\co 
(\varphi_{A}\ot H)\co (H\ot c_{H,A})\co (\delta_{H}\ot A)\co (\lambda_{H}\ot \lambda_{A})\co c_{A,H},$$
respectively, is a Hopf quasigroup if and only if (\ref{cesp1}) and 
\begin{equation}
\label{65}
(\phi_{A}\ot H)\co (A\ot ((\lambda_{H}\ot H)\co c_{H,H}\co \delta_{H}))=(\phi_{A}\ot H)\co (A\ot ((\lambda_{H}\ot H)\co \delta_{H}))
\end{equation}
hold.

Note that the product defined in (\ref{63}) is 
$$\mu_{A\circledast H}=(\mu_{A}\ot \mu_{H})\co (A\ot \Gamma \ot H),$$
where 
$$\Gamma=(\phi_{A}\ot H)\co (A\ot ((\lambda_{H}\ot H)\co c_{H,H}\co \delta_{H}))\co 
(\varphi_{A}\ot H)\co (H\ot c_{H,A})\co (\delta_{H}\ot A),$$
and 
$$\lambda_{A\circledast H}=\Gamma \co (\lambda_{H}\ot \lambda_{A})\co c_{A,H}.$$

It is easy to show that, if $(A, \varphi_{A}, \phi_{A})$  satisfies (\ref{61}),  (\ref{cesp2}),  (\ref{62}),  (\ref{cesp1}) and (\ref{65}),  we obtain that  the pairs $(A, \varphi_{A})$ and   $(A, \hat{\varphi}_{A}=\phi_{A}\co (A\ot \lambda_{H})\co c_{H,A})$ are  left $H$-quasimodule magmas-comonoids and satisfy the identities (\ref{cesp1}),  (\ref{cesp2}) and  (\ref{cesp4}). Moreover, under this conditions, 
$$
\Gamma=(\varphi_{A}\ot H)\co (H\ot T)\co (\delta_{H}\ot A), 
$$
where $T$ is defined as in (\ref{R2}). 

Therefore, by Example \ref{5}, the product defined in \cite{FT1} (see also  \cite[Theorem 3.9]{FT}) is  induced by an $a$-comonoidal distributive law of $H$ over $A$.
}
\end{example}

\section{The wreath product magma}

Let $A$ and $H$ be Hopf quasigroups. In this section, for an $a$-comonoidal distributive law of $H$ over $A$, $\Psi:H\ot A\rightarrow A\ot H$, we prove that $A\ot H$ with the  weak wreath product magma
\begin{equation}
\label{wpa1}
\mu_{A\ot_{\Psi}H}= (\mu_{A}\ot \mu_{H})\co (A\ot \Psi\ot H)
\end{equation}
becomes a Hopf quasigroup.

\begin{theorem}
\label{prin} Let $A$ and $H$ be Hopf quasigroups. Let $\Psi:H\ot A\rightarrow A\ot H$ be an $a$-comonoidal distributive law of $H$ over $A$. Then the  wreath product $A\ot_{\Psi} H$ built on $A\ot H$ with tensor coproduct and unit, where the product and antipode  are defined by (\ref{wpa1}) and 
\begin{equation}
\label{64}
\lambda_{A\ot_{\Psi} H}=\Psi\co (\lambda_{H}\ot \lambda_{A})\co c_{A,H},
\end{equation}
respectively, is a Hopf quasigroup. 
\end{theorem}

\begin{proof}
First note that, for $\eta_{A\ot_{\Psi}H}=\eta_{A}\ot \eta_{H}$, the identity
$$\mu_{A\ot_{\Psi}H}\co (\eta_{A\ot_{\Psi}H}\ot A\ot H)=id_{A\ot H}=\mu_{A\ot_{\Psi}H}\co (A\ot H\ot \eta_{A\ot_{\Psi}H})$$
follows from (\ref{dl3}), (\ref{dl4}) and the properties of the units. Then, the triple $(A\ot H, \eta_{A\ot_{\Psi}H}, \mu_{A\ot_{\Psi}H})$ is a unital magma.  On the other hand, $(A\ot H, \varepsilon_{A\ot_{\Psi}H}=\varepsilon_{A}\ot \varepsilon_{H}, \delta_{A\ot_{\Psi}H}=\delta_{A\ot H})$ is a comonoid, $\varepsilon_{A\ot_{\Psi}H}\co \eta_{A\ot_{\Psi}H}=id_{A\ot H}$ follows by (\ref{eta-eps}) for $H$ and $A$, $\varepsilon_{A\ot_{\Psi}H}\co \mu_{A\ot_{\Psi}H}=\varepsilon_{A\ot_{\Psi}H}\ot \varepsilon_{A\ot_{\Psi}H}$ follows by 
(\ref{mu-eps}) for $A$ and $H$ and by (\ref{cdl2}). Also,  $\delta_{A\ot_{\Psi}H}\co \eta_{A\ot_{\Psi}H}=\eta_{A\ot_{\Psi}H}\ot \eta_{A\ot_{\Psi}H}$ follows  by (\ref{delta-eta}) for $A$ and $H$ and by the naturality of $c$. 

On the other hand, 

\begin{itemize}
\item[ ]$\hspace{0.38cm}  (\mu_{A\ot_{\Psi}H}\ot \mu_{A\ot_{\Psi}H})\co \delta_{A\ot_{\Psi}H\ot A\ot_{\Psi}H}$

\item[ ]$= (\mu_{A}\ot \mu_{H\ot A}\ot \mu_{H})\co (A\ot A\ot c_{A,H}\ot c_{A,H}\ot H\ot H)\co (A\ot c_{A,A}\ot H\ot A\ot c_{H,H}\ot H)$
\item[ ]$\hspace{0.38cm}\co (\delta_{A}\ot ((\Psi\ot \Psi)\co \delta_{H\ot A})\ot \delta_{H})$ {\scriptsize ({\blue by naturality of $c$})} 

\item[ ]$= (\mu_{A}\ot \mu_{H\ot A}\ot \mu_{H})\co (A\ot A\ot c_{A,H}\ot c_{A,H}\ot H\ot H)\co (A\ot c_{A,A}\ot H\ot A\ot c_{H,H}\ot H)$
\item[ ]$\hspace{0.38cm}\co (\delta_{A}\ot ( \delta_{H\ot A}\co \Psi)\ot \delta_{H})${\scriptsize ({\blue by (\ref{cdl1})})} 

\item[ ]$= (A\ot c_{A,H}\ot H) \co (((\mu_{A}\ot \mu_{A})\co \delta_{A\ot A})\ot ((\mu_{H}\ot \mu_{H})\co \delta_{H\ot H}))\co (A\ot \Psi\ot H)$ {\scriptsize ({\blue by  naturality}} 
\item[ ]$\hspace{0.38cm}${\scriptsize {\blue  of $c$})}

\item[ ]$=\delta_{A\ot_{\Psi}H}\co \mu_{A\ot_{\Psi}H}$  {\scriptsize  ({\blue by (\ref{delta-mu}) for $A$ and $H$}).}

\end{itemize}

Moreover, if the antipode is defined as in (\ref{64}), we have 

\begin{itemize}
\item[ ]$\hspace{0.38cm} \mu_{A\ot_{\Psi}H}\circ (\lambda_{A\ot_{\Psi}H}\ot \mu_{A\ot_{\Psi}H})\circ (\delta_{A\ot_{\Psi}H}\ot A\ot H)$

\item[ ]$=(A\ot \mu_{H})\co (((\mu_{A}\ot H)\co (A\ot \Psi)\co (\Psi\ot A)\co (\lambda_{H}\ot \lambda_{A}\ot A))\ot H)\co (c_{A,H}\ot \mu_{A}\ot \mu_{H})$
\item[ ]$\hspace{0.38cm}\co (A\ot c_{A,H}\ot \Psi\ot H)\co (\delta_{A}\ot \delta_{H}\ot A\ot H)$ {\scriptsize ({\blue by definition})} 

\item[ ]$=(A\ot \mu_{H})\co ((\Psi\co (H\ot \mu_{A})\co (\lambda_{H}\ot \lambda_{A}\ot A))\ot H)\co (c_{A,H}\ot \mu_{A}\ot \mu_{H})\co (A\ot c_{A,H}\ot \Psi\ot H) $
\item[ ]$\hspace{0.38cm}\co (\delta_{A}\ot \delta_{H}\ot A\ot H)$ {\scriptsize ({\blue by (\ref{dl1})})} 

\item[ ]$=(A\ot \mu_{H})\co ((\Psi\co (\lambda_{H}\ot ( \mu_{A}\circ (\lambda_{A}\ot \mu_{A})\circ (\delta_{A}\ot A))))\ot \mu_{H})\co (c_{A,H}\ot \Psi\ot H)\co (A\ot \delta_{H}\ot A\ot H) $
\item[ ]$\hspace{0.38cm}$  {\scriptsize ({\blue  by naturality of $c$})} 

\item[ ]$=\varepsilon_{A}\ot ((A\ot \mu_{H})\co (\Psi\ot \mu_{H})\co (H\ot \Psi\ot H)\co (((\lambda_{H}\ot H)\co \delta_{H})\ot A\ot H)) $ {\scriptsize ({\blue by (\ref{lH}) and }} 
\item[ ]$\hspace{0.38cm}$ {\scriptsize {\blue naturality of $c$})}

\item[ ]$=\varepsilon_{A\ot_{\Psi} H}\ot A\ot H$  {\scriptsize  ({\blue by (\ref{adl1})}), }

\item[ ]

\end{itemize}

\begin{itemize}
\item[ ]$\hspace{0.38cm}\mu_{A\ot_{\Psi} H}\circ (A\ot H\ot \mu_{A\ot_{\Psi} H})\circ (A\ot H\ot \lambda_{A\ot_{\Psi} H}\ot A\ot H)\circ (\delta_{A\ot_{\Psi} H}\ot A\ot H) $

\item[ ]$=(\mu_{A}\ot \mu_{H})\co (A\ot \Psi\ot \mu_{H})\co (A\ot H\ot ((\mu_{A}\ot H)\co (A\ot \Psi)\co (\Psi\ot A)\co (\lambda_{H}\ot \lambda_{A}\ot A))\ot H) $
\item[ ]$\hspace{0.38cm} \co (A\ot H\ot c_{A,H}\ot A\ot H)\co (\delta_{A\ot H}\ot A\ot H)$ {\scriptsize ({\blue by definition})} 

\item[ ]$=(\mu_{A}\ot \mu_{H})\co (A\ot \Psi\ot \mu_{H})\co (A\ot H\ot (\Psi\co (H\ot \mu_{A})\co (\lambda_{H}\ot \lambda_{A}\ot A))\ot H)$
\item[ ]$\hspace{0.38cm}\co (A\ot H\ot c_{A,H}\ot A\ot H) \co (\delta_{A\ot H}\ot A\ot H)$ {\scriptsize ({\blue by (\ref{dl1})})} 

\item[ ]$= (\mu_{A}\ot H)\co (A\ot ( (A\ot \mu_{H})\co (\Psi\ot \mu_{H})\co (H\ot \Psi\ot H)\co (((H\ot \lambda_{H})\co \delta_{H})\ot A\ot H)))$
\item[ ]$\hspace{0.38cm}\co (A\ot H\ot (\mu_{A}\co (\lambda_{A}\ot A))\ot H)\co (A\ot c_{A,H}\ot A\ot H)\co (\delta_{A}\ot H\ot A\ot H)$ {\scriptsize ({\blue by  naturality}} 
\item[ ]$\hspace{0.38cm}$  {\scriptsize {\blue of $c$})}

\item[ ]$=((\mu_A\circ (A\ot \mu_A)\circ (A\ot \lambda_A\ot A)\circ (\delta_A\ot A))\ot H)\co (A\ot \varepsilon_{H}\ot A\ot H)$ {\scriptsize ({\blue by (\ref{adl2}) and  naturality}} 
\item[ ]$\hspace{0.38cm}$ {\scriptsize {\blue of $c$})}

\item[ ]$=\varepsilon_{A\ot_{\Psi} H}\ot A\ot H$  {\scriptsize  ({\blue by (\ref{lH}) for $A$}), }

\item[ ]

\end{itemize} 

\begin{itemize}
\item[ ]$\hspace{0.38cm}\mu_{A\ot_{\Psi} H}\circ (\mu_{A\ot_{\Psi} H}\ot A\ot  H)\circ (A\ot H\ot \lambda_{A\ot_{\Psi} H}\ot A\ot H)\circ (A\ot H\ot \delta_{A\ot_{\Psi} H}) $

\item[ ]$=(\mu_{A}\ot \mu_{H})\co (\mu_{A}\ot \Psi\ot H)\co (A\ot ((A\ot \mu_{H})\co ( \Psi\ot H)\co (H\ot \Psi)\co (H\ot \lambda_{H}\ot \lambda_{A}))\ot A\ot H)$
\item[ ]$\hspace{0.38cm}\co (A\ot H\ot c_{A,H}\ot A\ot H)\co (A\ot H\ot \delta_{A\ot H})$ {\scriptsize ({\blue by definition})} 

\item[ ]$= (\mu_{A}\ot \mu_{H})\co (\mu_{A}\ot \Psi\ot H)\co (A\ot (\Psi\co (\mu_{H}\ot A)\co (H\ot \lambda_{H}\ot \lambda_{A}))\ot A\ot H)$
\item[ ]$\hspace{0.38cm}\co (A\ot H\ot c_{A,H}\ot A\ot H)\co (A\ot H\ot \delta_{A\ot H})  ${\scriptsize ({\blue by (\ref{dl2})})} 

\item[ ]$= (A\ot \mu_{H})\co (((\mu_{A}\ot H)\co (\mu_{A}\ot \Psi)\co (A\ot \Psi\ot A)\co (A\ot H\ot (( \lambda_{A}\ot A)\co \delta_{A})))\ot H)$
\item[ ]$\hspace{0.38cm}\co (A\ot (\mu_{H}\co (H\ot \lambda_{H}))\ot A\ot H)\co (A\ot H\ot c_{A,H}\ot H)\co (A\ot H\ot A\ot \delta_{H})$ {\scriptsize ({\blue by naturality}}
\item[ ]$\hspace{0.38cm}${\scriptsize {\blue of $c$})}

\item[ ]$= (A\ot (\mu_H\circ (\mu_H\ot H)\circ (H\ot \lambda_H\ot H)\circ (H\ot \delta_H)))\co (A\ot H\ot \varepsilon_{A}\ot H)$ {\scriptsize ({\blue by (\ref{adl3}) and}} 
\item[ ]$\hspace{0.38cm}$ {\scriptsize {\blue naturality of $c$})} 

\item[ ]$= A\ot H\ot \varepsilon_{A\ot_{\Psi} H}$  {\scriptsize  ({\blue by (\ref{rH}) for $H$})}

\end{itemize} 

and 

\begin{itemize}
\item[ ]$\hspace{0.38cm} \mu_{A\ot_{\Psi} H}\circ(\mu_{A\ot_{\Psi} H}\ot \lambda_{A\ot_{\Psi} H})\circ (A\ot H\ot \delta_{A\ot_{\Psi} H})$

\item[ ]$= (\mu_{A}\ot H)\co (A\ot ((A\ot \mu_{H})\co ( \Psi\ot H)\co (H\ot \Psi)\co (H\ot \lambda_{H}\ot \lambda_{A})))\co (\mu_{A}\ot \mu_{H}\ot c_{A,H})$
\item[ ]$\hspace{0.38cm}\co (A\ot \Psi\ot c_{A,H}\ot H)\co (A\ot H\ot \delta_{A}\ot \delta_{H})$ {\scriptsize ({\blue by definition})} 

\item[ ]$=(\mu_{A}\ot H)\co (A\ot (\Psi\co (\mu_{H}\ot A)\co (H\ot \lambda_{H}\ot \lambda_{A}))\co (\mu_{A}\ot \mu_{H}\ot c_{A,H})\co (A\ot \Psi\ot c_{A,H}\ot H)$
\item[ ]$\hspace{0.38cm} \co (A\ot H\ot \delta_{A}\ot \delta_{H})  ${\scriptsize ({\blue by (\ref{dl2})})} 

\item[ ]$= (\mu_{A}\ot H)\co (A\ot \Psi)\co (\mu_{A}\ot (\mu_H\circ(\mu_H\ot \lambda_H)\circ (H\ot \delta_H))\ot \lambda_{A})\co (A\ot \Psi\ot c_{A,H})\co (A\ot H\ot \delta_{A}\ot H)$ 
\item[ ]$\hspace{0.38cm}$ {\scriptsize ({\blue by  naturality of $c$})}

\item[ ]$= ((\mu_{A}\ot H)\co (\mu_{A}\ot \Psi)\co (A\ot \Psi\ot A)\co (A\ot H\ot (( A\ot \lambda_{A})\co \delta_{A})))\ot \varepsilon_{H}$ {\scriptsize ({\blue by (\ref{rH}) for $H$ and }} 
\item[ ]$\hspace{0.38cm}$ {\scriptsize {\blue  naturality of $c$})}

\item[ ]$=A\ot H\ot \varepsilon_{A\ot_{\Psi} H}$  {\scriptsize  ({\blue by (\ref{adl4}) for $A$}).}

\item[ ]

\end{itemize} 

Therefore, (\ref{lH}) and (\ref{rH}) hold for $A\otimes_{\Psi}H$ and the proof is complete.

\end{proof}

\begin{example}
{\rm By Theorem \ref{prin} and the examples of the previous section, we obtain that the double cross product $A\bowtie H$ of two Hopf quasigroups $A$ and $H$, introduced in \cite{our2}, is an example of wreath product Hopf quasigroup. As a consequence, also is the cross product $A\bowtie_{\tau} H$ associated to a skew pairing (see \cite{our2}). Moreover, the Hopf quasigroups defined by the twisted double method in \cite{FT} are examples of wreath product Hopf quasigroups. Finally, the smash products of Hopf quasigroups introduced in \cite{BrzezJiao1} are examples of wreath product Hopf quasigroups as well as the twisted smash products defined in \cite{FT1}. 
}
\end{example}

\section*{Acknowledgements}
The  author was supported by  Ministerio de Ciencia e Innovaci\'on of Spain. Agencia Estatal de Investigaci\'on. Uni\'on Europea - Fondo Europeo de Desarrollo Regional (FEDER). Grant PID2020-115155GB-I00: Homolog\'{\i}a, homotop\'{\i}a e invariantes categ\'oricos en grupos y \'algebras no asociativas.


\begin{thebibliography}{99}

\bibitem{our1} Alonso Álvarez J. N., Fernández Vilaboa J. M., González Rodríguez R.,  Soneira Calvo, C., Projections  and Yetter-Drinfel'd modules over Hopf (co)quasigroups,  J. Algebra  443 (2015) 153-199.


\bibitem{our3}Alonso Álvarez, J.N., Fernández Vilaboa, J.M.,  González Rodríguez, R. y Soneira Calvo, C., Cleft comodules over  Hopf quasigroups, Communications in Contemporary Mathematics  17  (2015) 1550007.

\bibitem{our2} Alonso Álvarez J. N., Fernández Vilaboa J. M., González Rodríguez R., Multiplication alteration by two-cocycles. The  nonassociative version,  Bull. Malays. Math. Sci. Soc. 43 (2020),  3557-3615. 

\bibitem{our4} Alonso Álvarez J. N., Fernández Vilaboa J. M., González Rodríguez R., Quasigroupoids and weak Hopf quasigroups,   J. Algebra 568 (2021) 408-436.


\bibitem{Barr} Barr, M., Composite cotriples and derived functors, in: Sem. Triples Categor. Homology Theory, Springer
LN Math. 80 (1969), 336-356.

\bibitem{Beck} Beck, J., Distributive laws, in: Seminar on Triples and Categorical Homology Theory, B. Eckmann (ed.),
Springer LNM 80  (1969), 119-140.


\bibitem{BG1} B\"{o}hm, G., Gómez-Torrecillas, J., On the double crossed product of weak Hopf algebras,  AMS Contemp. Math. 585 (2013) 153-174.

\bibitem{BG2} B\"{o}hm, G., Gómez-Torrecillas, J., Bilinear factorization of algebras,  Bull. Belgian Math. Soc. 20 (2013) 221-244.

\bibitem{Bruck} Bruck R. H., Contributions to the theory of loops, Trans. Amer. Math. Soc.  60 (1946) 245-354.

\bibitem{BrzezJiao1} Brzezi\'nski T., Jiao Z. M., Actions of Hopf quasigroups, Comm. Algebra 40 (2012) 681-696.

\bibitem{BrzezJiao2} Brzezi\'nski T., Jiao Z. M., R-smash products of Hopf quasigroups, Arab. J. Math.  1 (2012) 39-46.

\bibitem{Chein} Chein O., Moufang loops of small order I,  Trans. Amer. Math. Soc.  188 (1974) 31-51.

\bibitem{FT1} Fang X., Wang, S. H.,  Twisted smash products for Hopf quasigroups, J. Southeast Univ.  27 (2011) 343-346.

\bibitem{FT} Fang X., Torrecillas B., Twisted smash products and L-R smash products for biquasimodule Hopf quasigroups,  Comm. Algebra 42 (2014) 4204-4234.

\bibitem{Christian} Kassel C., Quantum Groups,  GTM 155, Springer-Verlag, New York, 1995.

\bibitem{Majidesfera} Klim J., Majid S.; Hopf quasigroups and the algebraic 7-sphere, {\it J. Algebra } \textbf{323} (2010), 3067-3110.

\bibitem{SM2} Im, B., Nowak, A. W., Smith, J. D. H., Algebraic properties of quantum quasigroups,  J. Pure Appl. Algebra 225 (2021) 106539.

\bibitem{LP}   L\'opez L\'opez M. P., Villanueva Novoa E., The antipode and the (co)invariants of a finite Hopf (co)quasigroup,   Appl. Cat. Struct. 21 (2013) 237-247.

\bibitem{MAJDCP} Majid, S., Foundations of Quantum Group Theory, Cambridge University Press, Cambridge, 1995.

\bibitem{Mac} Mac Lane S., Categories for the Working Mathematician, GTM 5, Springer-Verlag, New-York, 1971.

\bibitem{PIS} P\'erez-Izquierdo J. M., Shestakov I. P., An envelope for Malcev algebras,  J. Algebra  272 (2004) 379-393.

\bibitem{PI07} P\'erez-Izquierdo J. M., Algebras, hyperalgebras, non-associative bialgebras and loops, Adv. Math.  208 (2007) 834-876.

\bibitem{SM1} Smith, J. D. H., Quantum quasigroups and loops, J. Algebra  456 (2016) 46-75.

\bibitem{RS} Street, R., Weak distributive laws, Theory Appl. Categ. 22 (2009) 313-320.




\end{thebibliography}
\end{document}